\documentclass[journal]{IEEEtran}
\usepackage{cite}
\usepackage{epstopdf}
\usepackage{subfigure}
\usepackage{bm}
\usepackage{amsmath}
\usepackage{cases} 
\usepackage{amssymb}
\usepackage{pifont}
\usepackage{xtab}
\usepackage{mathrsfs}
\usepackage{subeqnarray}
\usepackage{cases}
\usepackage{algorithm}
\usepackage{algcompatible}
\usepackage{booktabs}
\usepackage{multirow}
\usepackage{enumitem}
\usepackage{array}
\usepackage{diagbox}
\usepackage{algorithm}
\usepackage{algcompatible}
\usepackage[overload]{empheq}
\allowdisplaybreaks

\algblockdefx{FORALLP}{ENDFAP}[1]%
{\textbf{for all }#1 \textbf{do in parallel}}%
{\textbf{end for}}

\ifCLASSINFOpdf

\else

\usepackage[dvips]{graphicx}
\graphicspath{{../eps/}}
\DeclareGraphicsExtensions{.eps}
\fi

\IEEEoverridecommandlockouts

\newcolumntype{L}[1]{>{\raggedright\let\newline\\\arraybackslash\hspace{0pt}}m{#1}}
\newcolumntype{C}[1]{>{\centering\let\newline\\\arraybackslash\hspace{0pt}}m{#1}}
\newcolumntype{R}[1]{>{\raggedleft\let\newline\\\arraybackslash\hspace{0pt}}m{#1}}

\newcommand{\norm}[1]{\lVert#1\rVert}
\newcommand{\abs}[1]{\lvert#1\rvert}

\newtheorem {definition} {\textbf{Definition}}

\newtheorem {lemma}{\textbf{Lemma}}
\newtheorem {remark}{\textbf{Remark}}

\allowdisplaybreaks

\begin{document}
	
	\title{Distributionally Robust Chance-Constrained Approximate AC-OPF with Wasserstein Metric}
	
	\author{
		\IEEEauthorblockN{Chao Duan,~\IEEEmembership{Student Member,~IEEE}, Wanliang Fang, Lin Jiang,~\IEEEmembership{Member,~IEEE}, Li Yao,~\IEEEmembership{Student Member,~IEEE}, Jun Liu,~\IEEEmembership{Member,~IEEE}}
		% <-this % stops a space            
		\vspace{-1.0cm}
		
		\thanks{Manuscript received April 28, 2017; revised September 7, 2017 and December 27; accepted February 3, 2018. This work was supported in part by National Key Research and Development Program of China (2016YFB0901900), Physical Sciences Research Council (EP/L001004/1), and Key Research and Development Program of Shaanxi (2017ZDCXL-GY-02-03). \emph{(Corresponding author: Lin Jiang)}}
		
		\thanks{C. Duan, W. Fang, Li Yao and J. Liu are with the Department of Electrical Engineering, Xi'an Jiaotong University, Xi'an 710049, China. C. Duan is also with the Department of Electrical Engineering and Electronics, University of Liverpool, Liverpool L69 3GJ, U.K. (e-mail:duanchao@stu.xjtu.edu.cn; eewlfang@mail.xjtu.edu.cn; yaoli3021@sina.com; eeliujun@mail.xjtu.edu.cn)}
		
		\thanks{ L. Jiang is with Department of Electrical Engineering and Electronics, University of Liverpool, Liverpool L69 3GJ, U.K. (e-mail: ljiang@liverpool.ac.uk)}
		
	}
	
	\maketitle

	\begin{abstract}
		Chance constrained optimal power flow (OPF)  has been recognized as a promising framework to manage the risk from variable renewable energy (VRE). In presence of VRE uncertainties, this paper discusses a distributionally robust chance constrained approximate AC-OPF. The power flow model employed in the proposed OPF formulation combines an exact AC power flow model at the nominal operation point and an approximate linear power flow model to reflect the system response under uncertainties. The ambiguity set employed in the distributionally robust formulation is the Wasserstein ball centered at the empirical distribution. The proposed OPF model minimizes the expectation of the quadratic cost function w.r.t. the worst-case probability distribution and guarantees the chance constraints satisfied for any distribution in the ambiguity set. The whole method is data-driven in the sense that the ambiguity set is constructed from historical data without any presumption on the type of the probability distribution, and more data leads to smaller ambiguity set and less conservative strategy. Moreover, special problem structures of the proposed problem formulation are exploited to develop an efficient and scalable solution approach. Case studies are carried out on IEEE 14 and 118 bus systems to show the accuracy and necessity of the approximate AC model and the attractive features of the distributionally robust optimization approach compared with other methods to deal with uncertainties.

	\end{abstract}
	
	
	\IEEEpeerreviewmaketitle
	
	\vspace{-0.3cm}
	
	\section*{Notation}
	In this paper, boldface lower-case letter $\bm{x}$ represents a real vector and its $i$-th element is denoted by $x_{i}$. Boldface upper-case letter $\bm{A}$ represents a matrix with its $(i,j)^{\text{th}}$ element denoted by $A_{ij}$. Random vectors are written as boldface lower-case letter with tildes, i.e. $\tilde{\bm{y}}$. Given an index set $\mathcal{S}$, the subvector of $\bm{x}$ indexed by $\mathcal{S}$ is denoted by $\bm{x}_{\mathcal{S}}$. The imaginary unit is denoted by $j$. We use parentheses to construct vectors from comma separated lists as $(\bm{x}_,\bm{y},\bm{z})=[\bm{x}^{T},\bm{y}^{T},\bm{z}^{T}]^{T}$. In addition, the following special symbols are used in our problem formulation and derivation:    
	\addcontentsline{toc}{section}{Nomenclature}
	\begin{IEEEdescription}[\IEEEusemathlabelsep\IEEEsetlabelwidth{$V_1,V_2,V_3$}]
		\item[$n_b,n_l$] Number of buses and lines.
		\item[$\mathcal{R}$, $\mathcal{S}$, $\mathcal{L}$] Index sets for reference bus, PV buses and PQ buses, respectively.
		\item[$\bm{\theta},\bm{v}$] $n_b\times 1$ vectors of nominal bus voltage angles and magnitudes.
		\item[$\bm{f}$] $n_l\times 1$ vector of nominal line MW flow.
		\item[$\bm{Y},\bm{G},\bm{B}$] $\bm{Y}=\bm{G}+j\bm{B}$ is the system admittance matrix.
		\item[$\bm{B}'$] The susceptance matrix without shunt elements.
		\item[$\bm{G}^{l}$, $\bm{B}^{l}$] $n_{l}\times n_{b}$ matrices with ${{G}}^{l}_{ki}=-{{G}}^{l}_{kj}={G}_{ij}$, ${{B}}^{l}_{ki}=-{{B}}^{l}_{kj}={B}_{ij}$, and other elements being zeros.
		\item[$\bm{p},\bm{q}$] $n_b\times 1$ vectors of nominal active and reactive power injection.
		\item[$\bm{p}^{g},\bm{q}^{g}$] $n_b\times 1$ vectors of nominal active and reactive power injection from generators. $p^g_i=q^g_i=0$ if no generator at bus $i$.
		\item[$\bm{p}^{w},\bm{q}^{w}$] $n_b\times 1$ vectors of nominal active and reactive power injection from VRE.
		\item[$\bm{p}^{l},\bm{q}^{l}$] $n_b\times 1$ vectors of nominal active and reactive power consumption of load.
		\item[$\overline{\bm{p}}^g,\underline{\bm{p}}^g$] $n_b\times 1$ vectors of upper and lower output power limits of generators. $\overline{{p}}_i^g=\underline{{p}}_i^g=0$ if no generator at bus $i$.
		\item[$\bm{\alpha}$] $n_b\times 1$ vectors of AGC participation factors of generating units. $\alpha_i=0$ if no generator at bus $i$.
		\item[$\overline{\bm{r}},\underline{\bm{r}}$] $n_b\times 1$ vectors of upward and downward regulating reserves of generating units. $\overline{{r}}_i=\underline{{r}}_i=0$ if no generator at bus $i$.
		\item[$\overline{\bm{c}},\underline{\bm{c}}$] $n_b\times 1$ vectors of upward and downward regulating reserves prices.
		\item[$\tilde{\bm{\xi}}$] $n_b\times 1$ random vectors representing forecasting errors of VRE generation.
		\item[$\mathbb{P}$] A probability distribution (measure).
		\item[$\mathbb{E}_{\mathbb{P}}$] Expectation with probability distribution $\mathbb{P}$.
		\item[$\mathcal{P}(\Xi)$] Set of all distributions with support ${\Xi}$.
		\item[$(x)^{+}$] $\text{max}\{ x,0  \}$.
		\item[$\bm{A}\circ\bm{B}$] Hadamard product of matrix $\bm{A}$ and $\bm{B}$, i.e. $C_{ij}=A_{ij}B_{ij}$ if $\bm{C}=\bm{A}\circ\bm{B}$.
	\end{IEEEdescription}

	\section{Introduction}
	
	Large-scale VRE sources have been integrated into modern power systems with ever-increasing penetration. The uncertainties brought by VREs are no longer negligible and pose considerable risk to power system security. Risk management has been identified as one of the major challenges of integrating high penetrations of renewable energy \cite{challenge}. Therefore, developing optimization models and solution approaches with risk consideration is crucial for the reliable and economical operation of power systems. In recent years, we have seen extensive literature dealing with optimization of power system operation under uncertainties.
	
	Stochastic programming (SP) \cite{raey,5967923}, robust optimization (RO) \cite{6575173,lorca2015adaptive,wei2015robust,6948280} and distributionally robust optimization (DRO) \cite{7283666,7457327,7323859,7274485,6937216,7345592,7332992,7478165,8013071} have been employed to tackle uncertainties in power system operation. SP assumes the uncertainties follow a known probability distribution and transforms the stochastic problems into deterministic ones either by sampling or by analytical reformulation.  On the contrary, RO does not require any probabilistic information of the uncertainties. A deterministic uncertainty set is constructed to include all possible realizations of the random variable. RO seeks strategies that perform best w.r.t. the worst-case realization in the uncertainty set.  In fact, the uncertainties do obey some underlying probability distribution, but this distribution is not known as a priori and only some sample data is available. In other words, the underlying true distribution is ambiguous to the decision makers. To remedy the SP's specificity and RO's ignorance of the probabilistic information, DRO assumes that the true distribution lies in an ambiguity set and immunizes the operation strategies against all distributions in the ambiguity set. The most popular ambiguity set is moment-based, i.e. the set of all probability distributions with given mean and covariance \cite{7323859,7274485,6937216,7345592,7478165}. However, only the first two moments do not constitute a detailed characterization of the true distribution. Especially when we have a large amount of data at hand, much more probabilistic information can be extracted and exploited not just the first two moments. Intuitively, the more data is available, the more we know about the true distribution. Therefore, a desirable ambiguity set should be made smaller by incorporating more data. Moment-based ambiguity sets, obviously, do not possess such feature.
	
	In the framework of SP or DRO, the chance constraints, which require the security constraints hold with a specified probability level, have also been introduced to power system operation problems, especially the OPF problem \cite{bienstock2014chance,7332992,roald2015security,7478165,xie2016distributionally}. In \cite{bienstock2014chance,7332992,roald2015security}, by assuming the uncertainties follow Gaussian or Student's $t$ distribution, the chance constraints are cast into second-order cone programming (SOCP) constraints, and the uncertainties of the mean and covariance are further considered to achieve distributional robustness. In \cite{7478165,xie2016distributionally}, no assumption about the distribution type is used. The chance constraints are required to be satisfied for any probability distribution in the moment-based ambiguity set. All above existing works employ the DC power flow model without consideration for the voltage magnitudes and reactive power, which may lead to unpractical and unsafe operation strategy. Operated under the strategy obtained with DC power flow model, the system could be exposed to great risk of under- or over-voltage at some critical buses, and the generators would be forced to violate the excitation limits. Therefore, extending the chance-constrained OPF framework to include AC power flow model is a crucial step toward a better operation strategy under uncertainties.
	
	In this paper, we propose a Wasserstein-metric-based distributionally robust chance-constrained approximate AC-OPF. The contributions are as follows.   
	\begin{enumerate}[leftmargin=1\parindent]
		\item In the OPF formulation, to overcome the inaccuracy and V/Q unconcern of DC power flow model, we develop a hierarchical AC power flow model which is a combination of the exact nonlinear AC model at the nominal operation point and the approximate linear model developed in \cite{7782382} to represent the system response under uncertainties. It largely inherits the accuracy of the full AC model and also maintains the tractability of linear power flow model for use in stochastic optimization.   	
		
		\item We apply the recent results \cite{Esfahani2015Data,gao2016distributionally} of data-driven DRO with Wasserstein metric to the approximate AC-OPF model. To our knowledge, this is the first time the Wasserstein-metric-based ambiguity set is introduced to the OPF problems, which provides a promising alternative for the moment-based ambiguity set. The method is data-driven without any presumption on the probability distribution of the uncertainties. With only a limited number of data, the chance constraints w.r.t. the underlying true probability distribution can be robustly guaranteed. The more historical data is available, the less conservative the solution is.

		\item Note that the naive applications of the DRO approach discussed in \cite{Esfahani2015Data,gao2016distributionally} would lead to very unscalable implementations, i.e. the computational burden grows heavily with the number of data and dimension of uncertainties. To overcome such drawbacks, beyond a direct application of the general theories, the specific problem structures of the proposed OPF formulation are fully exploited to develop a scalable and efficient solution method.
	\end{enumerate}

	\section{Problem Formulation}
	
	\subsection{Power Flow and Its Control under Unceratainties}
\begin{figure}[!ht]\centering
		\includegraphics[width=3.2in]{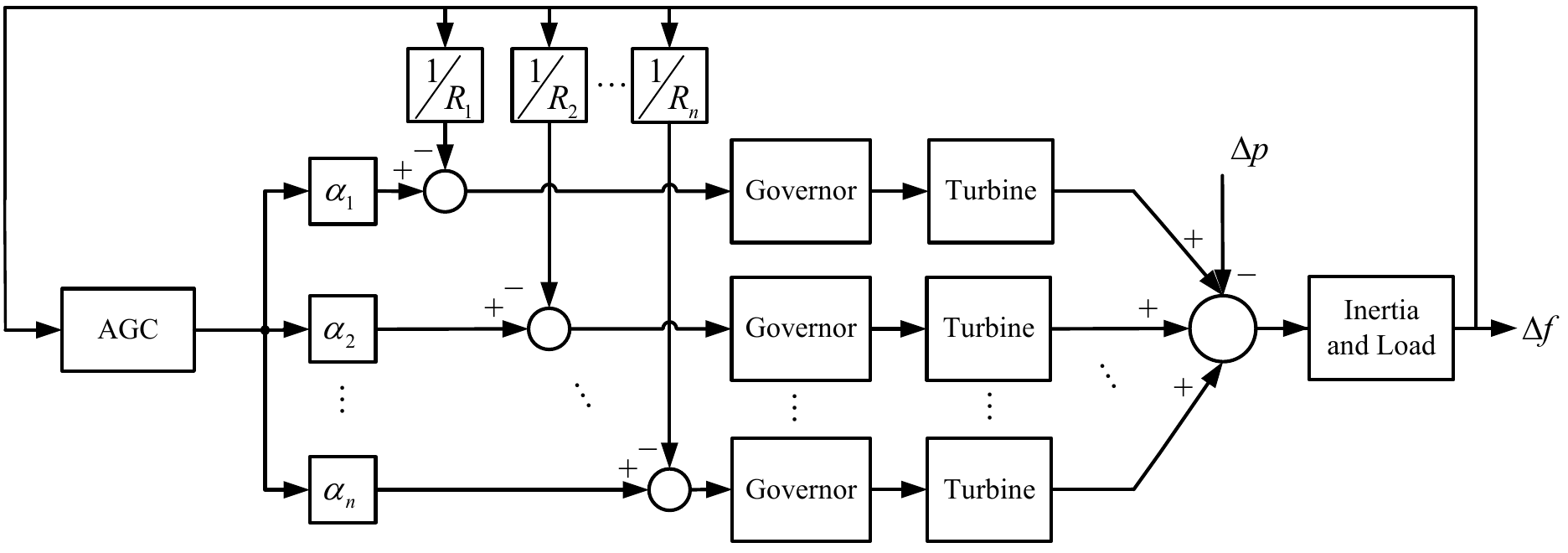}
		\caption{Automatic Generation Control} \label{AGC}
	\end{figure}

\begin{figure}[!ht]\centering
		\includegraphics[width=3.0in]{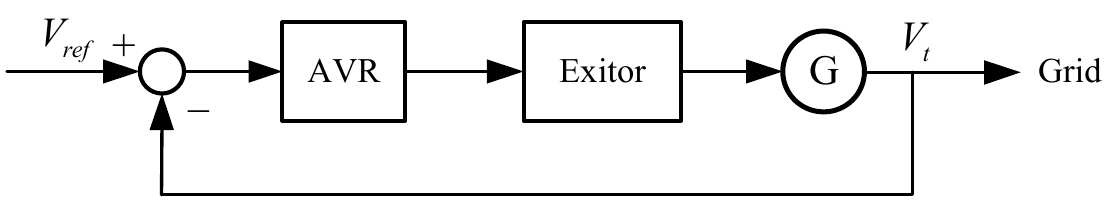}
		\caption{Automatic Voltage Regulation} \label{AVR}
	\end{figure}

	In a VRE integrated power system, the active and reactive power balance is governed by the following equations:
	\begin{subequations}\label{PF}
		\begin{align}[left =  \empheqlbrace\,]    
		&    \bm{p}^{g}+\bm{p}^{w}-\bm{p}^{l}=\mathscr{P}(\bm{\theta},\bm{v})\\
		&    \bm{q}^{g}+\bm{q}^{w}-\bm{q}^{l}=\mathscr{Q}(\bm{\theta},\bm{v})
		\end{align}
	\end{subequations}
	where 
	\begin{subequations}
		\begin{align}[left =  \empheqlbrace\,]    
		&    \mathscr{P}(\bm{\theta},\bm{v})=\left(\bm{G}\circ \bm{W} \circ\text{cos}\bm{\Theta}+\bm{B}\circ\bm{W} \circ\text{sin}\bm{\Theta}\right)\bm{1}\\
		&    \mathscr{Q}(\bm{\theta},\bm{v})=\left(-\bm{B}\circ \bm{W} \circ\text{cos}\bm{\Theta}+\bm{G}\circ\bm{W} \circ\text{sin}\bm{\Theta}\right)\bm{1}
		\end{align}
	\end{subequations}
	with $\bm{W}=\bm{v}\bm{v}^{\top}$ and $\bm{\Theta}=\bm{\theta}\bm{1}^{\top}-\bm{1}\bm{\theta}^{\top}$. For any fixed amount of renewable energy generation $\bm{p}^{w}+j\bm{q}^{w}$ and load demand $\bm{p}^{l}+j\bm{q}^{l}$, the conventional OPF addresses the problem of designing the optimal operation point $(\bm{\theta},\bm{v},\bm{p}^{g},\bm{q}^{g})$ which minimizes the generation costs while satisfying the system balance and security constraints. With increasing level of VRE penetration, the power injections across the network are facing continuous fluctuations due to the variability and unceratainty of VRE. Under such circumstance, a deterministic optimal operation point is no longer enough to guide the system operation. Therefore, the control mechanism to couple with renewable and load unceratinties must be taken into account in the OPF formulation. The most widely used control systems in response to volatility are the automatic generation control (AGC) and automatic voltage regulation (AVR) shown in Fig. \ref{AGC} and Fig. \ref{AVR} \cite{kundur1994power}. The AGC is a centralized control scheme to maintain the real-time active power balance by distributing the system power imbalance to each generating unit according to the corresponding participation factor. The AVR, on the other hand, is a decentralized control scheme to keep fixed voltage magnitudes at the buses equipped with reactive power scources. The most classic AVRs are those implemented in the generator exicitation systems, shown in Fig. \ref{AVR}, which maintain the stator voltages at the set values. With the implementation of the AGC and AVRs, the trasition of system operation point when renewable generation deviating from forecasting values is completely determined by the power flow equations (\ref{PF}). The new operation point can be computed by solving the conventional power flow problem with one generator bus as the reference bus, other generator buses as PV buses and the rest buses as PQ buses. Therefore, the new OPF formulation thus targets at deciding the optimal nominal operation point $(\bm{\theta},\bm{v},\bm{p}^{g},\bm{q}^{g})$ as well as the participation factors $\bm{\alpha}$ which provide a statistically economic and reliable system performance under uncertainties. However, the system response under uncertainties as described above is an inexplicit function given by the nonlinear power flow equations (\ref{PF}), which poses great challenge for formulating a tractable stochastic optimization problem. To circumvent this difficulty, an explicit surrogate formula is needed to approximate the exact implicit system response.

	\subsection{Linear Power FLow Response Model under Uncertainties}
	Recent advances in linear power flow model have opened a way for developing simple explicit formula for system response under uncertainties. The following linear power flow (LPF) model was developed in \cite{7782382}:
	\begin{equation}\label{LPF}
	\begin{bmatrix}
	\bm{p}\\ \bm{q}
	\end{bmatrix}
	=-
	\begin{bmatrix}
	\bm{B}' & -\bm{G} \\ 
	\bm{G}  & \bm{B}
	\end{bmatrix}
	\begin{bmatrix}
	\bm{\theta}\\ \bm{v}
	\end{bmatrix}
	\end{equation}
	
	\begin{equation}\label{brflow}
	\bm{f}=\bm{G}^{l}\bm{v}-\bm{B}^{l}\bm{\theta}
	\end{equation}
	where $\bm{p}$ and $\bm{q}$ are the $n_{b}\times 1$ vectors of nodal active and reactive power injections. The reference bus and PV buses are equipped with excitors or controllable reactive power compensators to maintain the pre-scheduled voltage magnitudes, and the reference bus is with fixed phase angle, typically 0. The state variables are arranged in the following sequence: $\bm{\theta}=(\bm{\theta}_{\mathcal{R}},\bm{\theta}_{\mathcal{S}},\bm{\theta}_{\mathcal{L}})$ and $\bm{v}=(\bm{v}_{\mathcal{R}},\bm{v}_{\mathcal{S}},\bm{v}_{\mathcal{L}})$. Other variables and coefficient matrices are ordered accordingly. By partially inverting (\ref{LPF}) and taking a incremental form, we have
	\begin{equation}\label{Vl}
	\begin{bmatrix}
	\Delta\bm{\theta}_{\mathcal{S}}\\ \Delta\bm{\theta}_{\mathcal{L}} \\ \Delta\bm{v}_{\mathcal{L}}
	\end{bmatrix}
	=-\bm{N}^{-1}
	\begin{bmatrix}
	\Delta\bm{p}_{\mathcal{S}}\\ \Delta\bm{p}_{\mathcal{L}} \\ \Delta\bm{q}_{\mathcal{L}}
	\end{bmatrix}
	-\bm{N}^{-1}\bm{H}
	\begin{bmatrix}
	0\\ \Delta\bm{v}_{\mathcal{R}} \\ \Delta\bm{v}_{\mathcal{S}}
	\end{bmatrix}
	\end{equation}
	
	\begin{equation}\label{Qg}
	\begin{bmatrix}
	\Delta\bm{q}_{\mathcal{R}}\\ \Delta\bm{q}_{\mathcal{S}} 
	\end{bmatrix}
	=\bm{L}\bm{N}^{-1}
	\begin{bmatrix}
	\Delta\bm{p}_{\mathcal{S}}\\ \Delta\bm{p}_{\mathcal{L}} \\ \Delta\bm{q}_{\mathcal{L}}
	\end{bmatrix}
	+(\bm{L}\bm{N}^{-1}\bm{H}-\bm{M})
	\begin{bmatrix}
	0\\ \Delta\bm{v}_{\mathcal{R}} \\ \Delta\bm{v}_{\mathcal{S}}
	\end{bmatrix}
	\end{equation}
	
	where 
	\begin{equation}
	\bm{H}
	=
	\begin{bmatrix}
	\bm{B}'_{\mathcal{SR}}    & -\bm{G}_{\mathcal{SR}} & -\bm{G}_{\mathcal{{SS}}} \\ 
	\bm{B}'_{\mathcal{LR}}    & -\bm{G}_{\mathcal{LR}} & -\bm{G}_{\mathcal{LS}} \\ 
	\bm{G}_{\mathcal{LR}}    & \bm{B}_{\mathcal{LR}} & \bm{B}_{\mathcal{LS}}
	\end{bmatrix}
	\end{equation}
	
	\begin{equation}
	\bm{N}
	=
	\begin{bmatrix}
	\bm{B}'_{\mathcal{SS}}    & \bm{B}'_{\mathcal{SL}} & -\bm{G}_{\mathcal{SL}} \\ 
	\bm{B}'_{\mathcal{LS}}    & \bm{B}'_{\mathcal{LL}} & -\bm{G}_{\mathcal{LL}} \\ 
	\bm{G}_{\mathcal{LS}}    & \bm{G}_{\mathcal{LL}} & \bm{B}_{\mathcal{LL}}
	\end{bmatrix}
	\end{equation}
	
	\begin{equation}
	\bm{M}
	=
	\begin{bmatrix}
	\bm{G}_{\mathcal{RR}}    & \bm{B}_{\mathcal{RR}} & \bm{B}_{\mathcal{RS}} \\ 
	\bm{G}_{\mathcal{SR}}    & \bm{B}_{\mathcal{SR}} & \bm{B}_{\mathcal{SS}} 
	\end{bmatrix}
	\end{equation}
	
	\begin{equation}
	\bm{L}
	=
	\begin{bmatrix}
	\bm{G}_{\mathcal{RS}}    & \bm{G}_{\mathcal{RL}} & \bm{B}_{\mathcal{RL}} \\ 
	\bm{G}_{\mathcal{SS}}    & \bm{G}_{\mathcal{SL}} & \bm{B}_{\mathcal{SL}} 
	\end{bmatrix}.
	\end{equation}    
	Given a nominal operation point $\bm{x}=(\bm{\theta},\bm{v},\bm{p}^{g},\bm{q}^{g})$, equation (\ref{Vl}) and (\ref{Qg}) constitute the explicit formulas for calculating system response under VRE and load uncertainties. Due to the implementation of AVR, (\ref{Vl}) and (\ref{Qg}) are further simplified as $\Delta\bm{v}_{\mathcal{R}}=\Delta\bm{v}_{\mathcal{S}}=\bm{0}$. The method in this paper can be used to consider any form of power injection uncertainties. But for ease of notation and without loss of generality, in the description of our method, we only consider the VRE uncertainties which are characterized by the random forecasting errors $\tilde{\bm{\zeta}}$ of VRE active power generation. The VREs are assumed to maintain a fixed power factor $\text{cos} \phi$ at the points of connection. Under the regulation of AGC, the generating units respond affinely to the total forecasting errors. As a result, the incremental active and reactive power injections in equation (\ref{Vl}) and (\ref{Qg}) become random variables (with tildes) as follows:
	\begin{subequations}\label{dealwithrandom}
		\begin{align}[left =  \empheqlbrace\,]    
		&\Delta\tilde{\bm{p}}_{\mathcal{S}}=-(\bm{1}^{\top}\tilde{\bm{\zeta}})\bm{\alpha}_{\mathcal{S}}+\tilde{\bm{\zeta}}_{\mathcal{S}}\\
		&\Delta\tilde{\bm{p}}_{\mathcal{L}}=\tilde{\bm{\zeta}}_{\mathcal{L}}\\
		&\Delta\tilde{\bm{q}}_{\mathcal{L}}=\sigma {\bm{\zeta}}_{\mathcal{L}}
		\end{align}
	\end{subequations}
	where $\sigma = \text{sin} \phi / \text{cos} \phi $ with $\text{cos} \phi$ being the mandated power factor of wind farms. 
	
	In summary, submitting (\ref{dealwithrandom}) into (\ref{Vl})(\ref{Qg}) along with (\ref{brflow}) and re-arranging the equations in an orgranized form, we obtain the following expression for system responses:
	\begin{subequations}\label{compactresponse}
		\begin{align}[left = \hspace{-0.2cm} \empheqlbrace\,]    
		& \tilde{\bm{v}}_{\mathcal{L}}=(\bm{1}^{\top}\tilde{\bm{\zeta}})\bm{A}^{v}\bm{\alpha}+\bm{B}^{v}\tilde{\bm{\zeta}}+\bm{v}_{\mathcal{L}}      
		\\
		& \tilde{\bm{q}}_{\mathcal{R}\cup \mathcal{S}}=(\bm{1}^{\top}\tilde{\bm{\zeta}})\bm{A}^{q}\bm{\alpha}+\bm{B}^{q}\tilde{\bm{\zeta}}+\bm{q}_{\mathcal{R}\cup \mathcal{S}}           
		\\
		&\tilde{\bm{f}}=(\bm{1}^{\top}\tilde{\bm{\zeta}})\bm{A}^{f}\bm{\alpha}+\bm{B}^{f}\tilde{\bm{\zeta}}+\bm{f}
		\end{align}
	\end{subequations}
	where $\bm{A}^{i},\bm{B}^{i}$, $i=v,q,f$ are all constant matrices decided by network parameters and nominal line MW flow $\bm{f}$ is given by (\ref{brflow}). To facilitate further discussion, we  rewrite (\ref{compactresponse}) in a component-wise form as follows:
	\begin{subequations}\label{compactresponse2}
		\begin{align}[left =  \empheqlbrace\,]    
		& \tilde{{v}}_{i}=    {{v}}_{i}(\bm{x},\bm{1}^{\top}\tilde{\bm{\zeta}},\bm{B}^{v}_{i:}\tilde{\bm{\zeta}}), i \in \mathcal{L}   
		\\
		& \tilde{{q}}_{j}=    {{q}}_{i}(\bm{x},\bm{1}^{\top}\tilde{\bm{\zeta}},\bm{B}^{q}_{j:}\tilde{\bm{\zeta}}), j \in \mathcal{R}  \cup   \mathcal{S}
		\\
		&\tilde{{f}}_{k}={{f}}_{k}(\bm{x},\bm{1}^{\top}\tilde{\bm{\zeta}},\bm{B}^{f}_{k:}\tilde{\bm{\zeta}}), k=1,2,\cdots,n_l
		\end{align}
	\end{subequations}
	where $\bm{B}^{v}_{i:}$ denotes the $i$th row of matrix $\bm{B}^{v}$ and this notation extends to other matrices as well; the explicit formulas for ${{v}}_{i}(\cdot)$, ${{q}}_{j}(\cdot)$ and ${{f}}_{k}(\cdot)$ can be easily seen from (\ref{compactresponse}).

	\subsection{Problem Formulation of Distributionally Robust Chance-constrained Approximate AC-OPF}
	
	Based on the discussion in the last subsection, we formulate the distributionally robust chance-constrained approximate AC-OPF as follows:
	\begin{equation}\label{WDROPF_obj}
	\underset{    
		\bm{x}
	}{\text{min}} \ \underset{\mathbb{P}\in \hat{\mathcal{P}}_{N}}{\text{sup}}
	\mathbb{E}_{\mathbb{P}} \left\{ \sum_{i =1}^{n_b} f_{i}\left({{p}^g_{i}}-(\bm{1}^{\top}\tilde{\bm{\zeta}}){\alpha}_{i}\right) + \overline{\bm{c}}^{\top}\overline{\bm{r}}+\underline{\bm{c}}^{\top}\underline{\bm{r}} \right\}
	\end{equation}        \vspace{-0.0cm}
	\begin{subequations}\label{WDROPF_con}\hspace{-2cm}
		\begin{align}[left = \text{s.t.} \empheqlbrace\,]    
		&    \bm{p}^{g}+\bm{p}^{w}-\bm{p}^{l}=\mathscr{P}(\bm{\theta},\bm{v})\label{con0}\\
		&    \bm{q}^{g}+\bm{q}^{w}-\bm{q}^{l}=\mathscr{Q}(\bm{\theta},\bm{v})\label{con1}\\
		& \underline{\bm{v}} \leq \bm{v} \leq \overline{\bm{v}}\label{conv}\\
		&\bm{1}^{\top}\bm{\alpha}=1, \ \bm{\alpha} \geq \bm{0} \label{con2}\\
		&\bm{p}^g-\underline{\bm{r}}\geq \underline{\bm{p}}^g,\ \underline{\bm{r}}\geq \bm{0} \label{con3}\\
		& \bm{p}^g+\overline{\bm{r}}\leq \overline{\bm{p}}^g, \ \overline{\bm{r}}\geq \bm{0}\label{con4}\\
		&\underset{\mathbb{P}\in \hat{\mathcal{P}}_{N}}{\text{inf}} \mathbb{P}\left\{-\underline{\bm{r}} \leq -(\bm{1}^{\top}\tilde{\bm{\zeta}})\bm{\alpha} \leq \label{con5}  \overline{\bm{r}}\right\}\geq 1-\rho_1 \\
		&  \forall i \in \mathcal{L} ,\ j \in \mathcal{R}  \cup   \mathcal{S}, \   k=1,2,\cdots,n_l : \nonumber\\
		&\underset{\mathbb{P}\in \hat{\mathcal{P}}_{N}}{\text{inf}} \mathbb{P}\left\{ \underline{v}_i \leq \tilde{v}_i\leq \overline{v}_i  \right\}\geq 1-\rho_2 \label{con6}\\
		&\underset{\mathbb{P}\in \hat{\mathcal{P}}_{N}}{\text{inf}} \mathbb{P}\left\{ \underline{q}_{j} \leq \tilde{q}_{j} \leq \overline{q}_{j}  \right\}\geq 1-\rho_3 \label{conq}\\
		&\underset{\mathbb{P}\in \hat{\mathcal{P}}_{N}}{\text{inf}} \mathbb{P}\left\{ \underline{f}_{k} \leq \tilde{f}_{k} \leq \overline{f}_{k} \right\}\geq 1-\rho_4 \label{con7}
		\end{align}
	\end{subequations}    
	where $\bm{x}=(\bm{v},\bm{\theta},\bm{p}^g,\bm{q}^g,\bm{\alpha},\overline{\bm{r}},\underline{\bm{r}})$; $\mathcal{P}_N$ is an ambiguity set of probability distributions constructed from historical data and its explicit expression will be discussed in the next section. The objective function (\ref{WDROPF_obj}) is the worst-case expectation of generation and reserve costs in which the cost function $f_{i}(\cdot), i=1,\cdots,n_g$ are convex quadratic functions of the form:
	\begin{equation}\label{costfun}
	f_{i}(x)= c_{i2}x^2+c_{i1}x+c_{i0}.
	\end{equation}
	with $c_{i2},c_{i1},c_{i0} \geq 0$.  Nolinear equality constraint (\ref{con0}) and (\ref{con1}) ensure the active and reactive power balance at the nominal opeartion point where the VRE generation exactly match the forecasting values. Constraint (\ref{conv}) sets the voltage magnitude limits for all buses at the nominal opeartion point and also for the reference and PV buses under perturbed conditions. Constraint (\ref{con2}) enforces the basic requirement for participation factors in the AGC system. Constraint (\ref{con3}) and (\ref{con4}) guarantee the reserve availability considering the output limits of generating units. The adequacy of downward and upward regulating reverse is ensured by distributionally robust chance constraint (\ref{con5}). Similarly, constraint (\ref{con6}) ensures the voltage quality for PQ buses, and constraint (\ref{conq}) safeguards the adequancy of reactive power for AVR at generator buses. Finally, constraint (\ref{con7}) guarantees the adequacy of transmission line capacity. 
	
	\begin{remark}
		The OPF formulation (\ref{WDROPF_obj})(\ref{WDROPF_con}) combines the exact nonlinear power flow equations (\ref{PF}) and the linear approximate system response model (\ref{compactresponse2}). On the one hand, this new problem formulation inherits the basic feature of conventional AC-OPF by using the exact nonlinear power flow to govern the system state at the nominal condition. On the other hand, the linear approximate power flow response model is employed to describe the system behavior when VRE generation deviates from forecasting values, which makes it possible to develop efficient techniques to deal with the stochastic formulation (\ref{WDROPF_obj})(\ref{con5})$\sim$(\ref{con7}). Of course, though more tractable, the linear power flow response model still brings inaccuracies in evaluating the operational costs, voltages, reactive power and line flows. The accuracy of this approach will be assessed by comparative numerical studies in section V.
	\end{remark}
	
	\begin{remark}
		Note that the dimension of random variable $\tilde{\bm{\zeta}}$ is equal to the number of VREs installed across the network. It could be prohibitively high for numerical computation in a real-world large-scale power system with distributed wind farms and PVs. Fortunately, the simple affine policy employed in the AGC scheme enables significant reduction in the dimension of the random variable considered in model (\ref{WDROPF_obj})(\ref{WDROPF_con}). Specifically, in the objective function (\ref{WDROPF_obj}) and constraint (\ref{con5}), we only need to consider the 1-dimensional random variable $\bm{1}^{\top}\tilde{\bm{\zeta}}$. For the constraint (\ref{con6}), it is only required to consider the 2-dimensional random variable $(\bm{1}^{\top}\tilde{\bm{\zeta}},\bm{B}^{v}_{i:}\tilde{\bm{\zeta}})$ in light of equation (\ref{compactresponse2}). Similarly, the 2-dimensional random variable $(\bm{1}^{\top}\tilde{\bm{\zeta}},\bm{B}^{f}_{i:}\tilde{\bm{\zeta}})$ needs to be taken into account in constraint (\ref{con7}). When the historical data of $\tilde{\bm{\zeta}}$ is available at hand, it is equivalent to have the historical data for the 1- or 2-dimensional random variables described above. In summary, no matter how many uncertain renewable sources are installed in the system, we only need to deal with 1- or 2-dimensional random variables which significantly improve the scalability of the proposed method.
	\end{remark}

\begin{remark}
In a conventional stochastic chance-constrained programming, the probability distribution $\mathbb{P}$ is assumed to be known apriori. However, in practice, the probability distribution of $\tilde{\xi}$ is unknown and only some historical data is available. Theoretically, the precise knowledge of the probability distribution cannot be obtained from finite data. Therefore, the probability distribution is ambiguous. Nonetheless, the historical data does provide us some reliable probabilistic information based on which we can construct an ambiguity set $\mathcal{P}_{N}$, i.e. a set of probability distributions consistent with the observed historical data. Therefore, the objective function (\ref{WDROPF_obj}) and the chance-constraints (\ref{con5})$\sim$(\ref{con7}) consider the worst-case distribution in the ambiguity set $\mathcal{P}_{N}$ to ensure the strategy performs well under the ambiguity of distribution.
\end{remark}

	\section{Ambiguity Set Using Wasserstein Metric}
	
	\subsection{Wasserstein Metric and Ambiguity Set}
	
	As shown in equation (\ref{compactresponse}), the system responses are the functions of some random variable $\tilde{\bm{\xi}}$. Whether evaluating the expected generation costs or assessing the reliability of system operation, we need to have a knowledge about the probability distribution $\mathbb{P}$ of random variable $\tilde{\bm{\xi}}$ under consideration. Unfortunately, in the real-life application, we only have a finite set of historical data and the precise characterization of $\mathbb{P}$ can never be extracted from finite samples. In other words, the underlying true distribution $\mathbb{P}$ is ambiguous. However, according to the historical sample set $\{\hat{\bm{\xi}}^{(1)},\hat{\bm{\xi}}^{(2)},\cdots,\hat{\bm{\xi}}^{(N)} \}$ at hand, we can construct an empirical distribution $\mathbb{\hat{P}}_{N}=\frac{1}{N}\sum_{k=1}^{N}\delta_{\hat{\bm{\xi}}^{(k)}}$ which acts as an estimation of the true distribution $\mathbb{P}$. Intuitively, $\mathbb{\hat{P}}_{N}$ converges to $\mathbb{P}$ as $N \rightarrow \infty$, i.e. the ``distance'' between $\mathbb{\hat{P}}_{N}$ and $\mathbb{P}$ becomes smaller when more data is available. One of the ``distances'' to establish the convergence of $\mathbb{\hat{P}}_{N}$ to $\mathbb{P}$ is the Wasserstein metric defined as follows:    
	\begin{definition}[Wasserstein Metric]
		For any probability distribution $\mathbb{Q}_{1},\mathbb{Q}_{2} \in \mathcal{P}(\Xi)$, the Wasserstein metric can be defined through  
		\begin{equation}
		\begin{aligned}
		W(\mathbb{Q}_{1},&\mathbb{Q}_{2} )=\underset{\Pi}{\text{inf}}\bigg\{ \int_{\Xi^{2}} \norm{\xi_{1}-\xi_{2}} \Pi (d\xi_2, d\xi_2):\\
		&
		\begin{aligned}
		&\Pi\ \text{is\ a\ joint\ distribution\ of}\ \xi_1\ \text{and}\ \xi_2\\ 
		&\text{with\ marginal distributions}\ 
		\mathbb{Q}_{1}\ \text{and}\ \mathbb{Q}_{2}
		\end{aligned}
		\bigg\}.
		\end{aligned} 
		\end{equation}
	\end{definition}
	where $\mathcal{P}(\Xi)$ denotes the set of all probability distributions with support $\Xi$; $\norm{\cdot}$ can be any norm in $\mathbb{R}^n$, and we use $l_{1}$ norm $\norm{\cdot}_{1}$ in this paper for its superior numerical tractability in DRO. Hence, we have $W(\mathbb{\hat{P}}_{N} , \mathbb{P})\leq \epsilon(N)$ where $\epsilon(\cdot)$ is some sample-dependent monotone function decreasing to $0$ as $N$ tends to infinity, and the explicit formula for $\epsilon(\cdot)$ will be discussed in the next subsection. Therefore, given a historical data set with $N$ samples, we know that the true distribution $\mathbb{P}$ belongs to the following set:
	\begin{equation}\label{amset}
	\hat{\mathcal{P}}_{N}=\left\{ \mathbb{P}\in \mathcal{P}(\Xi): W(\mathbb{P} ,\mathbb{\hat{P}}_{N}) \leq \epsilon (N)  \right\}
	\end{equation}    
	which is called a \emph{ambiguity set} of the underlying true distribution. Note that $\hat{\mathcal{P}}_{N}$ is the Wasserstein ball of radius $\epsilon (N)$ centered at the empirical distribution $\mathbb{\hat{P}}_{N}$. It represents the reliable information about the true distribution $\mathbb{P}$ observed from the $N$ historical samples at hand.

	\subsection{Selection of Wasserstein Radius $\epsilon (N)$}
	
	The radius of the Wasserstein ball in the ambiguity set (\ref{amset}) is crucial for the performance of the distributionally robust chance-constrained optimization problem (\ref{WDROPF_obj})(\ref{WDROPF_con}). Following the discussion in last subsection, $\epsilon (N)$ is a decreasing function, as small as possible, to bound $W(\mathbb{\hat{P}}_{N} , \mathbb{P})$ from above. One possible choice for $\epsilon (N)$ was given in \cite{Zhao2015Data} as follows:
	\begin{equation}\label{eps}
	\epsilon (N) = D \sqrt{\frac{2}{N} \text{log}\left(\frac{1}{1-\beta}\right)},
	\end{equation}
	where $\beta$ is a confidence level and $D$ is the diameter of the support of the random variable. However, in our numerical experience, we find this choice of $\epsilon (N)$ overly conservative for the DRO to gain noticable advantage over the conventional RO with reasonable amount of data.
	
	To improve the formula (\ref{eps}), let's review the derivation of (\ref{eps}) presented in \cite{Chaoyue2015Data,Zhao2015Data}. In the proof of proposition 3 in \cite{Chaoyue2015Data}, it is shown that for any $\delta >0$
	\begin{equation}\label{concentrate}
	\mathbb{P}[W(\mathbb{\hat{P}}_{N} , \mathbb{P}) \geq \epsilon] \leq \text{exp}\left(
	-N\underset{\mathbb{Q}:W(\mathbb{Q} , \mathbb{P}) \geq \epsilon}{\text{inf}}(
	H(\mathbb{Q}|\mathbb{P})-\delta
	)
	\right)
	\end{equation}
	where $H(\mathbb{Q}|\mathbb{P})$ is called the Kullback information of $\mathbb{Q}$ w.r.t. $\mathbb{P}$ (see \cite{MR2172583} for more detail). Corollary 2.4 in \cite{MR2172583} shows
	\begin{equation}\label{talgad}
	H(\mathbb{Q}|\mathbb{P}) \geq W(\mathbb{Q},\mathbb{P})^{2}/C^{2}, \ \forall \mathbb{Q} \in \mathcal{P}(\Xi)
	\end{equation}
	where 
	\begin{equation}\label{Cvalue}
	C=2\underset{\bm{\xi}_{0} \in \Xi, \alpha >0}{\text{inf}}\left(
	\frac{1}{2\alpha}(1+\text{ln}\mathbb{E}_{\mathbb{P}}[e^{\alpha \norm{\tilde{\bm{\xi}}-\bm{\xi}_{0}}_{1}^{2}}])
	\right)^{1/2}
	\end{equation}
	with a particular case $C=\sqrt{2}D$ (particular case 2.5 in \cite{MR2172583}). Plugging (\ref{talgad}) into RHS of (\ref{concentrate}) and taking $\delta \rightarrow 0$, we have
	\begin{equation}\label{concentrate2}
	\mathbb{P}[W(\mathbb{\hat{P}}_{N} , \mathbb{P}) \leq \epsilon]\geq 1-\text{exp}\left(
	-N \frac{\epsilon^{2}}{C^{2}}
	\right).
	\end{equation}
	Equating the RHS of (\ref{concentrate2}) to $\beta$ leads to
	\begin{equation}\label{eps2}
	\epsilon (N) = C \sqrt{\frac{1}{N} \text{log}\left(\frac{1}{1-\beta}\right)}.
	\end{equation}
	Note that if the special case $C=\sqrt{2}D$ is adopted, (\ref{eps2}) turns into (\ref{eps}). However, $\sqrt{2}D$ is much larger than (\ref{Cvalue}) in practice, which is one of the major sources of conservatism for (\ref{eps}). Instead of adopting the special case $C=\sqrt{2}D$, we directly use (\ref{Cvalue}) by estimating its value from data. From (\ref{Cvalue}), we have
	\begin{equation}
	\begin{aligned}
	C \leq & 2\underset{\alpha >0}{\text{inf}}\left(
	\frac{1}{2\alpha}(1+\text{ln}\mathbb{E}_{\mathbb{P}}[e^{\alpha \norm{\tilde{\bm{\xi}}-\hat{\bm{\mu}}}_{1}^{2}}])
	\right)^{1/2} \\
	\approx & 2\underset{\alpha >0}{\text{inf}}\left(
	\frac{1}{2\alpha}\left(1+\text{ln}(
	\frac{1}{N}\sum_{k=1}^{N}e^{\alpha \norm{\hat{\bm{\xi}}^{(k)}-\hat{\bm{\mu}}}_{1}^{2}})\right)
	\right)^{1/2} 
	\end{aligned}
	\end{equation}
	where $\hat{\bm{\mu}}$ denotes the sample mean and the minimization over $\alpha$ can be easily done by the bisection search method.

	\section{Solution Approach}
	
	To solve the distributionally robust voltage-concerned chance-constrained AC-OPF (\ref{WDROPF_obj})(\ref{WDROPF_con}), we need to evaluate the worst-case expected costs in (\ref{WDROPF_obj}) and reformulate the distributionally robust chance constraints (\ref{con5})$\sim$(\ref{con7}). The underlying engine for the reformulation is the strong duality result developed in \cite{gao2016distributionally}. We tailor the result to our problem setting as the following lemma.
	\begin{lemma}[\cite{gao2016distributionally}]
		Given a random variable $\tilde{\bm{\xi}} \in \mathbf{R}^{m}$ with closed and convex support $\Xi$, the Wasserstein ball $\mathcal{B}_{\epsilon}(\mathbb{\hat{P}}_{\xi})$ is constructed from sample set $\{\hat{\bm{\xi}}_{1},\hat{\bm{\xi}}_{2},\cdots,\hat{\bm{\xi}}_{N} \}$. If the loss function $l(\tilde{\bm{\xi}})$ is upper semi-continuous, the worst-case expectation
		\begin{equation}
		\begin{aligned}
		&\underset{\mathbb{P}\in \hat{\mathcal{P}}_{N}}{\text{sup}}\
		\mathbb{E}_{\mathbb{P}}\left\{
		l(\tilde{\bm{\xi}})
		\right\}\\
		=& \underset{\lambda \geq 0}{\text{inf}} \left\{
		\ \lambda \cdot \epsilon+\frac{1}{N}\sum_{k=1}^{N}\underset{\bm{\xi} \in \Xi}{\text{sup}} \left( l(\bm{\xi})-\lambda \norm{\bm{\xi}-\hat{\bm{\xi}}_{k}}_{1}  \right)\right\}.
		\end{aligned}
		\end{equation}    
	\end{lemma}
Note that a naive application of Lemma 1 will make the optimization problem scale grows with number of historical data, which is computationally undesirable. To develop a more efficient solution approach, using Lemma 1 while exploiting the problem structures, this section replaces the worst-case expected cost with its tight upper bound (Section-IV-A) and safely approximates the chance constraints by corresponding robust constraints (Section-IV-B). By doing so, the computational burden almost remains unchanged as more historical data is available at hand.
	
	In our implementation, we estimate the support $\Xi$ from the sample set as
	\begin{equation}
	\Xi =\{\bm{\xi} \in \mathbf{R}^{m} | -\sigma_{\text{max}}\bm{1}\leq \hat{\bm{\Sigma}}^{-1/2}(\bm{\xi}-\hat{\bm{\mu}}) \leq \sigma_{\text{max}}\bm{1}  \}
	\end{equation}
	where $\hat{\bm{\mu}}$ and $\hat{\bm{\Sigma}}$ are the sample mean and sample covariance, and $\sigma_{\text{max}}=10$ in our practice.

	\subsection{Evaluation of Worst-case Costs}
	
	By noticing (\ref{costfun}), the cost function inside the $\mathbb{E}_{\mathbb{P}}\{\cdot \}$ operator in (\ref{WDROPF_obj}) can be written as follows:
	\begin{equation}
	\eta(\bm{x},\tilde{\omega})=c_2\tilde{\omega}^2-c_1\tilde{\omega}+c_0
	\end{equation}
	where
	\begin{equation}\label{ccc}
	\left\{
	\begin{array}{l}
	c_2=\sum_{i=1}^{n_{g}}c_{i2}{\alpha}_i^2\\
	c_1=\sum_{i=1}^{n_{g}} 2c_{i2}{p}_{i}{\alpha}_i+c_{i1}{\alpha}_i\\
	c_0=\sum_{i=1}^{n_{g}} (c_{i2}{p}_{i}^{2}+c_{i1}{p}_{i}+c_{i0})+ \overline{\bm{c}}^{\top}\overline{\bm{r}}+\underline{\bm{c}}{\top}\underline{\bm{r}}.
	\end{array}\right.
	\end{equation}
	with $\tilde{\omega}=\bm{1}^{\top}\tilde{\bm{\zeta}}$. Since $c_{i2} \geq 0$, we have $c_2 \geq 0$. It follows that $\eta(\bm{x},\cdot)$ is a convex quadratic function. In addition, $\eta(\cdot,\tilde{\omega})$ is also a convex quadratic function because it is the sum of convex quadratic funtions $\{f_{i}\}_{i=1}^{n_g}$.
	
	Given sample set $\{\hat{\omega}_{1},\hat{\omega}_{2},\cdots,\hat{\omega}_{N} \}$ of random variable $\tilde{\omega}$ and its support $[\underline{\omega},\overline{\omega}]$, we can evaluate the worst-case costs in (\ref{WDROPF_obj}) using Lemma 1 as follows:
	\begin{subequations}
		\begin{align}
		& 
		\underset{\mathbb{P}\in \hat{\mathcal{P}}_{N}}{\text{sup}}\
		\mathbb{E}_{\mathbb{P}}\left\{
		\eta(\bm{x},\tilde{\omega})
		\right\}
		\\ \label{wc_1}
		=& \left\{
		\begin{array}{l}
		\underset{\lambda \geq 0,\bm{s} \in \mathbf{R}^{N}}{\text{inf}} 
		\ \lambda \cdot \epsilon+\frac{1}{N}\sum_{k=1}^{N}s_{k}\\
		\text{s.t.}\hspace{-0.1cm} \underset{\underline{\omega}\leq \omega \leq \overline{\omega}}{\text{sup}}\hspace{-0.2cm} \left( \eta(\bm{x},\omega)-\lambda \abs{\omega-\hat{\omega}_{k}}  \right) \leq s_{k},\forall k\leq N
		\end{array}\right.
		\\ \label{wc_2}
		=&  \left\{
		\begin{array}{l}
		\underset{\lambda \geq 0,\bm{s} \in \mathbf{R}^{N}}{\text{inf}} 
		\ \lambda \cdot \epsilon+\frac{1}{N}\sum_{k=1}^{N}s_{k}\\
		\text{s.t.}\  \eta(\bm{x},\underline{\omega})+\lambda ({\underline{\omega}-\hat{\omega}_{k}})  \leq s_{k},\forall k\leq N\\
		\quad \ \eta(\bm{x},\overline{\omega})-\lambda ({\overline{\omega}-\hat{\omega}_{k}})  \leq s_{k},\forall k\leq N\\
		\quad \ \eta(\bm{x},\hat{\omega}_{k})  \leq s_{k},\forall k\leq N
		\end{array}\right.
		\end{align}
	\end{subequations}
	where (\ref{wc_1}) follows from a direct application of Lemma 1 and the introduction of epigraphical auxiliary variables $s_{k}, k\leq N$; equality (\ref{wc_2}) comes from the observation that $\eta(\bm{x},\omega)-\lambda \abs{\omega-\hat{\omega}_{k}}$ is convex in $\omega$ in the intervals $[{\underline{\omega},\hat{\omega}_{k}}]$ and $[\hat{\omega}_{k},\overline{\omega}]$, repectively. Hence, the supremum in (\ref{wc_1}) can only be attained at the boundary of the intervals, i.e. $\underline{\omega}$, $\hat{\omega}_{k}$ or $\overline{\omega}$. Note that the worst-case expected costs have already been transformed into a tractable deterministic formulation (\ref{wc_2}) with linear objective function and convex quadratic constraints. However, problem (\ref{wc_2}) has a huge computational disadvantage: the numbers of auxiliary variables and quadratic constraints are proportional to the size of historical sample set. When a large data set is at hand, the computational burden could prevent us fully exploiting the historical data.
	
	To overcome this drawback, we replace (\ref{wc_2}) with a close upper approximation which is more scalable w.r.t. the size of sample set. Let $\eta ' (\bm{x},\tilde{\omega})=2c_2\tilde{\omega}-c_1$ be the derivative of $\eta  (\bm{x},\tilde{\omega})$ w.r.t. $\tilde{\omega}$, and take $\lambda = \text{max}\{\eta ' (\bm{x},\overline{\omega}),-\eta ' (\bm{x},\underline{\omega})  \}$ which is the smallest value of $\lambda$ in (\ref{wc_2}) such that 
	\begin{equation}\label{etaleq}
	\left\{
	\begin{array}{l}
	\eta(\bm{x},\underline{\omega})+\lambda ({\underline{\omega}-{\omega}}) \leq \eta(\bm{x},{\omega}),\ \forall {\omega}\in [\underline{\omega},\overline{\omega}]\\
	\eta(\bm{x},\overline{\omega})-\lambda ({\overline{\omega}-{\omega}})  \leq \eta(\bm{x},{\omega}),\ \forall {\omega}\in [\underline{\omega},\overline{\omega}].
	\end{array}\right.
	\end{equation}
	Thus,
	\begin{equation}\label{wc_3}
	(\ref{wc_2}) \leq \left\{
	\begin{array}{l}
	\underset{\lambda \in \mathbf{R}}{\text{inf}} 
	\ \lambda \cdot \epsilon+\frac{1}{N}\sum_{k=1}^{N}\eta(\bm{x},\hat{\omega}_{k})\\
	\text{s.t.}\  \eta'(\bm{x},\overline{\omega})  \leq \lambda\\
	\quad \  -\eta'(\bm{x},\underline{\omega})  \leq \lambda.
	\end{array}\right.
	\end{equation}
	Note that the numbers of decision variables and constraints in (\ref{wc_3}) stay unchanged when using larger sample set, therefore (\ref{wc_3}) is more computationally favorable than (\ref{wc_2}). In addition, the optimums of (\ref{wc_2}) and (\ref{wc_3}) have very small difference in practice and such difference diminishes as the sample size $N$ increases.
	\vspace{-0.1cm}

	\subsection{Reformulation of Distributionally Robust Chance Constraints}
	
	The chance constraints (\ref{con5})$\sim$(\ref{con7}) can be written in the general form:
	\begin{equation}\label{DRCC}
	\underset{\mathbb{P}\in \hat{\mathcal{P}}_{N}}{\text{inf}} \mathbb{P}[\bm{g}(\bm{x},\tilde{\bm{\xi}})\leq \bm{0}]\geq 1-\rho.
	\end{equation}
	where $\bm{g}$ depends linearly on decision variable $\bm{x}$ and random variable $\tilde{\bm{\xi}}$, respectively.
	The above constraint is in fact non-convex, so it is generally difficult to derive a tractable equivalent reformulation. In this paper, we develop a convex conservative approximation to (\ref{DRCC}) through the following strategy. We seek a deterministic uncertainty set $\mathcal{U}$ such that the robust constraint
	\begin{equation}\label{RC}
	\bm{g}(\bm{x},{\bm{\xi}})\leq \bm{0},\ \forall {\bm{\xi}} \in \mathcal{U}
	\end{equation}
	implies the distributionally robust chance constraint (\ref{DRCC}). In addition, the uncertainty set $\mathcal{U}$ should possess the following features: 1) obtaining $\mathcal{U}$ is computationally cheap; 2) robust constraint (\ref{RC}) is numerically tractable and efficient; 3) uncertainty set $\mathcal{U}$ is made as small as possible.
	
\begin{figure}[!ht]\centering
		\includegraphics[width=3.4in]{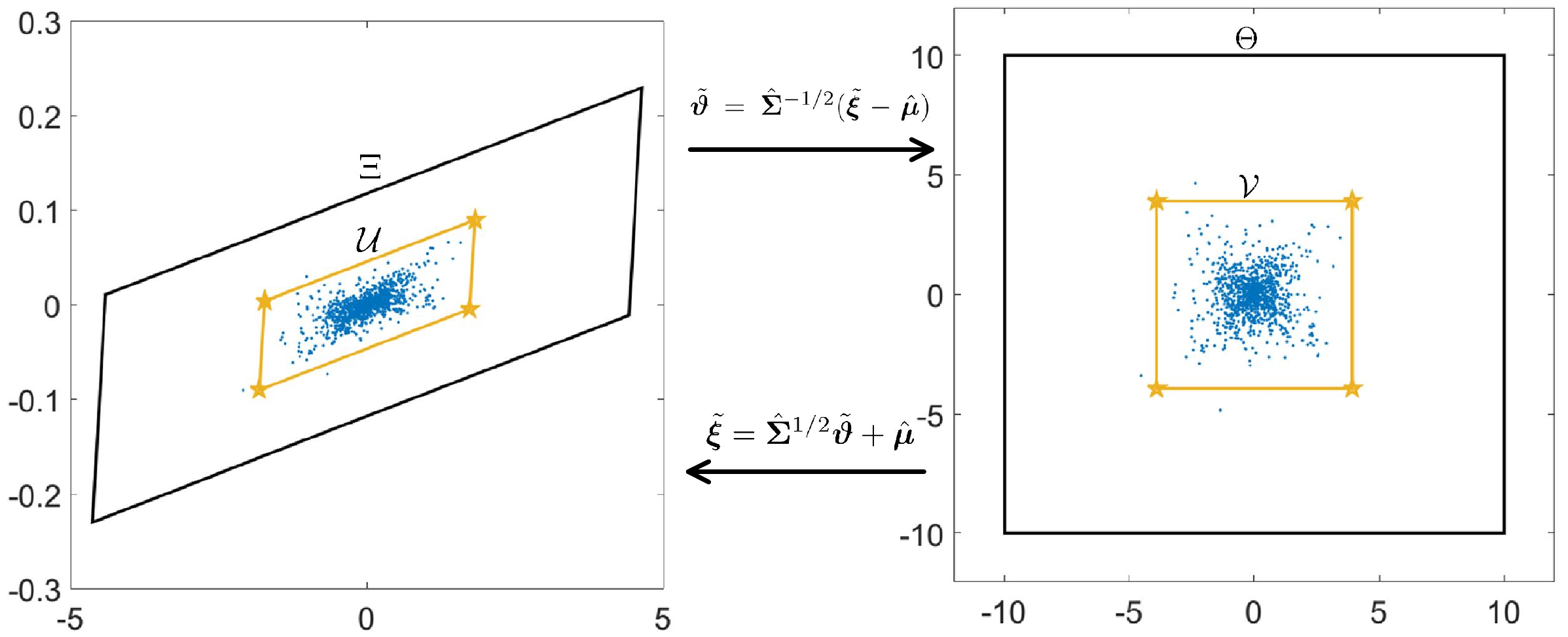}
		\caption{The procedure to determine the uncertainty set $\mathcal{U}$ in two dimension.} \label{uncertainset}
	\end{figure}

	Given a sample set     $\{\hat{\bm{\xi}}^{({1})},\hat{\bm{\xi}}^{({2})},\cdots,\hat{\bm{\xi}}^{({N})} \}$ of random variable $\tilde{\bm{\xi}} \in \mathbf{R}^{m}$ with the unknown underlying true distribution $\mathbb{P}$, we can compute the sample mean $\hat{\bm{\mu}}$ and the sample covariance $\hat{\bm{\Sigma}}$. Instead of directly working with random variable $\tilde{\bm{\xi}}$, we consider its standardized version $\tilde{\bm{\vartheta}}=\hat{\bm{\Sigma}}^{-1/2}(\tilde{\bm{\xi}}-\hat{\bm{\mu}})$ with sample set $\{\hat{\bm{\vartheta}}^{(k)}=\hat{\bm{\Sigma}}^{-1/2}(\hat{\bm{\xi}}^{(k)}-\hat{\bm{\mu}}) \}_{k=1,2,\cdots,N}$,  illustrated in Fig. \ref{uncertainset}. Obviously, random variable $\tilde{\bm{\vartheta}}$ has sample mean $\bm{0}$,  sample covariance $\mathbf{I}$ and support $\Theta = \{ -\sigma_{\text{max}}\bm{1} \leq \bm{\vartheta} \leq \sigma_{\text{max}}\bm{1} \}$. Let $\mathbb{Q}$ and $\mathbb{Q}_{N}$ denote the true distribution and empirical distribution of $\tilde{\bm{\vartheta}}$, respectively. Accordingly, the ambiguity set $\hat{\mathcal{Q}}_{N}$ can be constructed as in (\ref{amset}). We seek a set $\mathcal{V} \subseteq \mathbf{R}^{m}$ such that
	\begin{equation}
	\underset{\mathbb{Q}\in \hat{\mathcal{Q}}_{N}}{\text{sup}} \mathbb{Q}[\tilde{\bm{\vartheta}} \notin \mathcal{V}]\leq \rho.
	\end{equation}
	Then $\mathcal{U}=\hat{\bm{\Sigma}}^{1/2} \mathcal{V}+\hat{\bm{\mu}}$ would be a desired uncertainty set needed in (\ref{RC}). To develop an efficient method to find out such $\mathcal{V}$, considering the sample uncorrelatedness and equal variance among different components of $\tilde{\bm{\vartheta}}$, we restrict the structure of $\mathcal{V}$ as the hypercube
	\begin{equation}
	\mathcal{V}(\sigma )=\{\bm{\vartheta} \in \mathbf{R}^{m} | -\sigma \bm{1} <  \bm{\vartheta} <   \sigma \bm{1}  \}
	\end{equation}
	whose side length $\sigma$ needs to be minimized in order to reduce conservatism. This leads to the following optimization problem:
	\begin{equation}\label{minside}
	\underset{0\leq \sigma \leq \sigma_{\text{max}} }{\text{min}}\  \sigma \ \ \text{s.t.} \underset{\mathbb{Q}\in \hat{\mathcal{Q}}_{N}}{\text{sup}} \mathbb{Q}[\tilde{\bm{\vartheta}} \notin \mathcal{V}(\sigma)]\leq \rho.
	\end{equation}
	By leveraging the worst-case probability formulation in \cite{Esfahani2015Data} and the duality result in Lemma 1,
	problem (\ref{minside}) has a very simple deterministic reformulation which is a consequence of the following lemma.
	\begin{lemma}
		\begin{equation}\label{wcprob}
		\begin{aligned}
		&\underset{\mathbb{Q}\in \hat{\mathcal{Q}}_{N}}{\text{sup}} \mathbb{Q}[\tilde{\bm{\vartheta}} \notin \mathcal{V}(\sigma)]\\
		=&\underset{\lambda \geq 0}{\text{inf}} \left\{
		\  \lambda \cdot \epsilon+\frac{1}{N}\sum_{k=1}^{N}\left(1-\lambda (\sigma-\norm{\hat{\bm{\vartheta}}^{({k})}}_{\infty})^{+}  \right)^{+} \right\}.    
		\end{aligned}
		\end{equation}
	\end{lemma}
	where $(x)^{+}=\text{max}(x,0)$.
	\begin{IEEEproof}
		The proof is delayed to the appendix.
	\end{IEEEproof}
	Using Lemma 2, problem (\ref{minside}) is equivalent to 
	\begin{equation}\label{minside2}
	\underset{0 \leq \lambda ,0\leq \sigma \leq \sigma_{\text{max}}}{\text{min}}  \sigma \ \ \text{s.t.} \ \  h(\sigma,\lambda) \leq \rho.
	\end{equation}
	where $h(\sigma,\lambda) = \lambda\epsilon+\frac{1}{N}\sum_{k=1}^{N}\left(1-\lambda (\sigma-\norm{\hat{\bm{\vartheta}}^{({k})}}_{\infty})^{+}  \right)^{+}$.
	
	Note that the worst-case probability (\ref{wcprob}) is non-decreasing in $\sigma$, therefore problem (\ref{minside}), equivalently problem (\ref{minside2}), has a unique minimum. Although problem (\ref{minside2}) is non-smooth, it only involves two 1-dimensional decision variables. We can design a nested bisection search method to quickly locate the optimal solution. The method is summarized in Algorithm 1 in which the function $\texttt{bisearch}(f(\cdot),a,b)$ returns the minimum of $f(\cdot)$ in the interval $[a,b]$ by performing a bisection search. Further note that $h(\sigma,\lambda)$ is convex in $\lambda$ for fixed $\sigma$, so the bisection search in step 4 of Algorithm 1 is well-defined. Since Algrithm 1 only involves function evaluations, it solves problem (\ref{minside2}) very efficiently.

	\begin{algorithm}
		\caption{Nested Bisection Search}
		\begin{algorithmic}[1]
			\STATE Initialize $\underline{\sigma}=0$, $\overline{\sigma}=\sigma_{\text{max}}$;
			\WHILE{($\overline{\sigma}-\underline{\sigma}>10^{-4}$)}
			\STATE $\sigma = (\overline{\sigma}+\underline{\sigma})/2$; 
			\STATE $\gamma=$\texttt{bisearch}$(h(\sigma,\cdot)$,0,100);
			\IF{$\gamma > \rho$}
			\STATE $\underline{\sigma}=\sigma$; 
			\ELSE
			\STATE $\overline{\sigma}=\sigma$;
			\ENDIF
			\ENDWHILE
			\STATE Output $\sigma=\overline{\sigma}$.
		\end{algorithmic}
	\end{algorithm}
	
	After determining the optimal $\sigma$, the hypercube $\mathcal{V}(\sigma )$ can be expressed as the convex hull of its vertices. For 1-dimensional random variable, $\mathcal{V}(\sigma )=\text{conv}(\{-\sigma, \sigma  \})$, and for 2-dimensional random variable, $\mathcal{V}(\sigma )=\text{conv}(\{(\pm \sigma,\pm \sigma) \})$. In general, $\mathcal{V}(\sigma )=\text{conv}(\{ \bm{v}^{(1)},\cdots,\bm{v}^{(2^{m})} \})$. Hence the desired uncertainty set $\mathcal{U}$ takes the form
	\begin{equation}
	\mathcal{U}=\text{conv}(\{ \bm{u}^{(1)},\cdots,\bm{u}^{(2^{m})} \})
	\end{equation}
	where $\bm{u}^{(i)}=\hat{\bm{\Sigma}}^{1/2} \bm{v}^{(i)}+\hat{\bm{\mu}}$, $1 \leq i \leq 2^{m}$. Thus, noticing the linear dependence of $\bm{g}$ on $\bm{\xi}$, the robust constraint (\ref{RC}) is equivalent to
	\begin{equation}\label{LEC}
	\bm{g}(\bm{x},{\bm{u}}^{(i)})\leq \bm{0},\ 1 \leq i \leq 2^{m}
	\end{equation}
	which is a set of linear inequality constraints on decision variable $\bm{x}$. As a result, the distributionally robust chance constraint (\ref{DRCC}) can be safely approximated by the linear constraints (\ref{LEC}). Note that we only need to deal with the cases of $m=1$ and $m=2$ as commented in Remark 2, therefore constraint (\ref{LEC}) does not bring any scalability problems and can be handled very efficiently. The whole procedure is illustrated in Fig. \ref{uncertainset}.
	
	To sum up, replacing the worst-case expected costs in (\ref{WDROPF_obj}) with  (\ref{wc_3}), taking the place of chance constraints (\ref{con5})$\sim$(\ref{con7}) with their corresponding deterministic constraint (\ref{LEC}), we obtain a deterministic optimization problem with quadratic objective function (\ref{WDROPF_con}), nonlinear equality constraints (\ref{con0})(\ref{con1}), and linear constraints (\ref{conv})$\sim$(\ref{con7}). In fact, this problem structure can be seen as an extension of conventional AC-OPF by adding decision varibles $(\bm{\alpha},\overline{\bm{r}},\underline{\bm{r}})$ and a set of linear constraints. Therefore, the problem can be solved (at least to local optimality) by a mature interior point method (IPM) solver.

	\section{Numerical Results}
	This section reports numerical results. We will first discuss some implementation issues of the proposed method and the related benchmarking approaches. Then the case studies on IEEE 14 and 118 bus systems will be presented. The tests on 14-bus system emphasize on validity and necessity of the proposed approximate AC-OPF formulation. Whereas the studies on 118-bus system focus on the features of the DRO approach with the comparison to other methods to deal with uncertainties.
	
	\subsection{Implementation and Benchmarking}
	The MATPOWER \cite{zimmerman2011matpower} offers a very convenient way to implement our problem formulation. It employs an extensible OPF structure \cite{zimmerman2009matpower} to allow the user to modify and augment the problem formulation without rewriting the portions that are shared with the conventional AC-OPF. Our problem formulation fits well into this structure and the modification needed can be easily done with the help of YALMIP \cite{lofberg2004yalmip}. We choose KNITRO \cite{byrd2006knitro} as the default IPM solver. Further note that voltage magnitude (\ref{con6}), reactive power (\ref{conq}) and line flow (\ref{con7}) constraints are only active at a fraction of buses, generators and lines in practical systems. To improve the computational efficiency, we employ the successive constraint enforcement scheme \cite{6948280} which relaxes all the constraints (\ref{con6})$\sim$(\ref{con7}) initially and subsequently adds back those that are violated at the solution of the relaxed model.   
	
	In order to comprehensively assess the performance of the proposed Wasserstein-metric-based distributionally robust optimization (\textbf{WDRO}) method, three other methods are introduced for benchmarking. The first one is the robust optimization (\textbf{RO}) approach which requires the security constraints (\ref{con5})$\sim$(\ref{con7}) to be satisfied for all realizations in the support of the random variable. The second benchmarking approach is the moment-based distributionally robust optimization (\textbf{MDRO}). The ambiguity set for MDRO is the set of all probability distributions with given mean and covariance (usually sample mean and sample covariance). Then the chance constraints (\ref{con5})$\sim$(\ref{con7}) can be reformulated as SOCP constraints by leveraging the Chebyshev inequality $\mathbb{P}\{ \abs{\tilde{\xi}-\mu} \geq \sqrt{{1}/{\rho}}\sigma \}\leq \rho$. Moreover, the third approach is the Gauss-based stochastic programming (\textbf{GSP}) which presumes that the random variable follows Gauss distribution with given mean and covariance. Then the chance constraints (\ref{con5})$\sim$(\ref{con7}) can also be reformulated into SOCP constraints by using the inequality $\mathbb{P}\{ \abs{\tilde{\xi}-\mu} \geq \Phi^{-1}(1-\rho/2)\sigma \}\leq \rho$ where $\Phi$ is the cumulative distribution function of standard Gauss random variable. In all three benchmarking approaches described above, the expectation in the objective function (\ref{WDROPF_obj}) is estimated by the sample average of the cost function.    It seems that the MDRO and GSP cannot be integrated into the extensible OPF structure provided by MATPOWER due to the introduction of SOCP constraints. However, rather than directly working with the SOCP constraints, the cutting-plane approach illustrated in \cite{7332992} and \cite{bienstock2014chance}  iteratively employs a sequence of linear cuts to the SOCP constraints. In each iteration, only linear constraints are involved in enforcing the chance constraints (\ref{con5})$\sim$(\ref{con7}), therefore it fits well in the extensible OPF structure provided by MATPOWER. 
	
	To test and compare the performance of different methods to deal with uncertainties, an underlying random number generator (RNG) is employed to simulate the VRE forecasting errors. Following the suggestion in \cite{Lee2016Probabilistic}, we use Laplace distribution to generate ``realistic'' historical data for wind power forecasting errors with the typical standard error in U.S. reported in \cite{wind2012}. The parameters of RNG including the type of distribution is secret and only the data generated from the RNG is available to all the methods under test. From a finite number of data, each method constructs its own optimal operation strategy. Then the RNG is again used to generate a much larger set of data for Monte Carlo simulation (MCS) to assess the statistical performance of the operation strategy obtained by each method.
	
	\subsection{IEEE 14-Bus System}
\begin{figure}[!ht]\centering
		\includegraphics[width=3.0in]{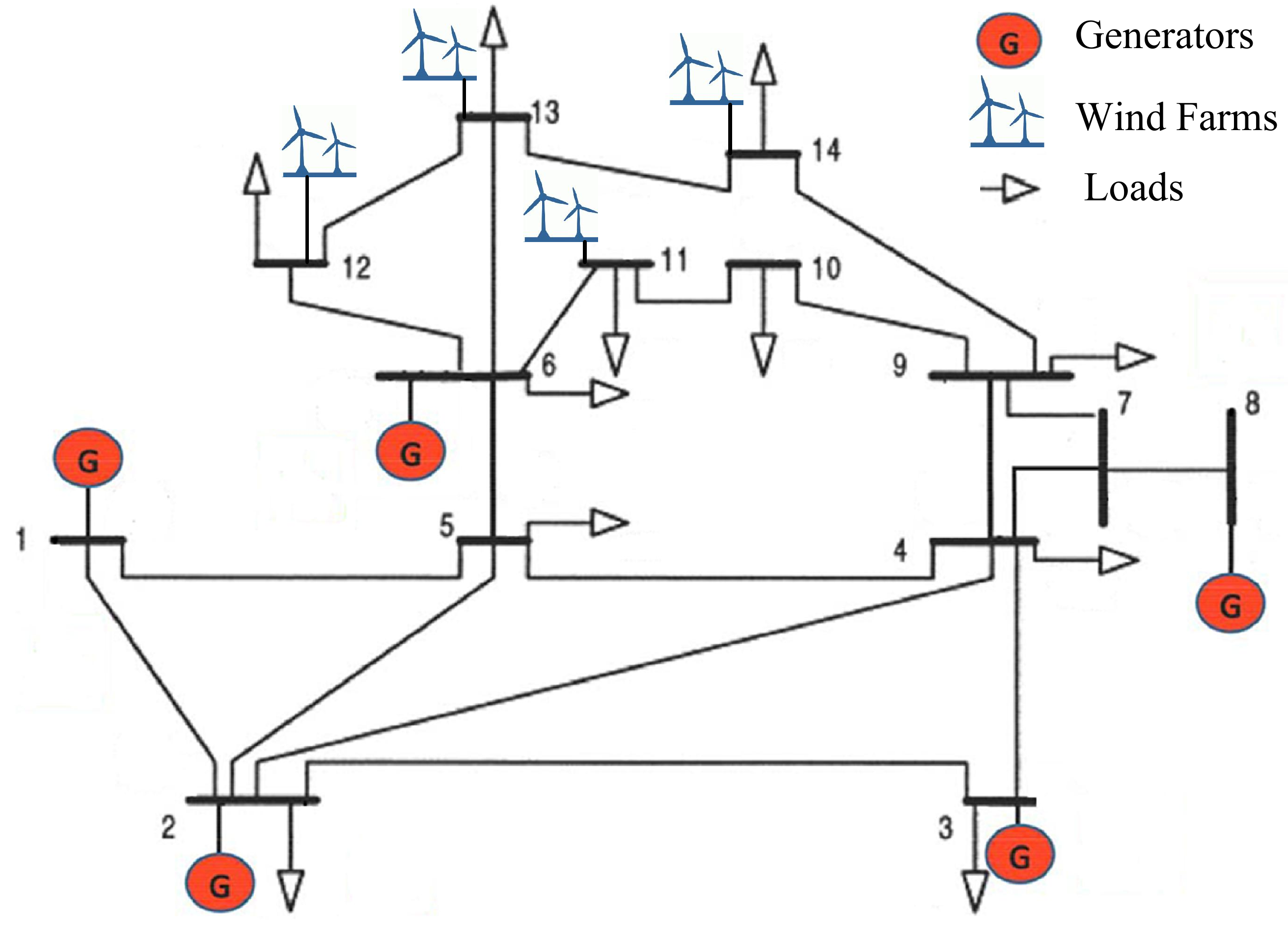}
		\caption{Diagram of the Modified IEEE 14-bus System} \label{ieee14sys}
	\end{figure}

	The diagram of the modified IEEE 14-bus system is shown in Fig. \ref{ieee14sys}. Four wind farms with each capacity of 36 MW are installed at bus 11, 12, 13 and 14. The forecasting values of wind power are $50\%$ of their capacity. The transmission capacity of each line is set to be 40 MW. The regulating reserve prices are assumed to be $50\%$ of the linear coefficients of the generator cost functions. The following group of tests is designed to demonstrate the features of the proposed approximate AC-OPF model.
	
\begin{figure}[!ht]\centering
		\includegraphics[width=1.7in]{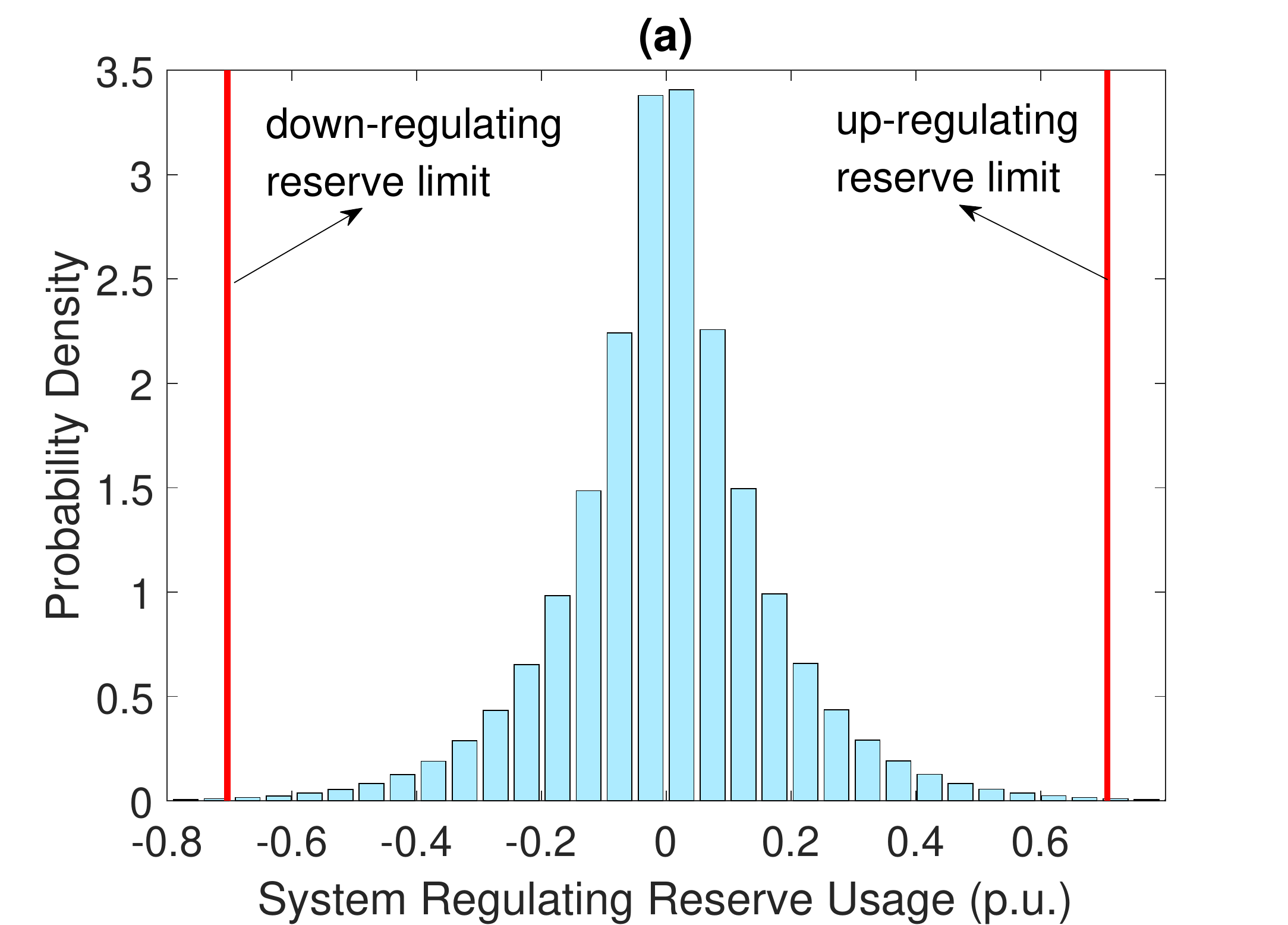}
		\includegraphics[width=1.7in]{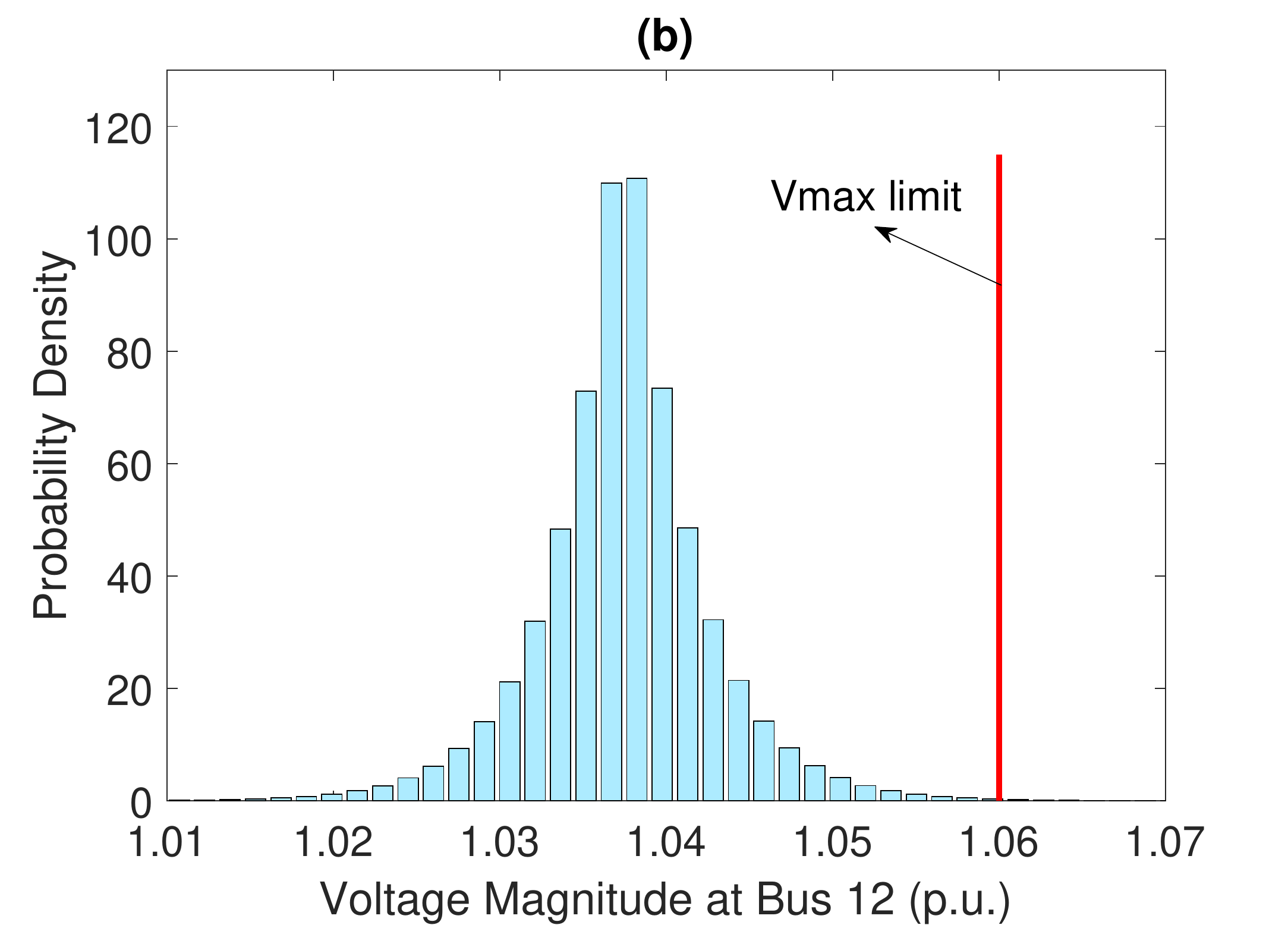}
		\includegraphics[width=1.7in]{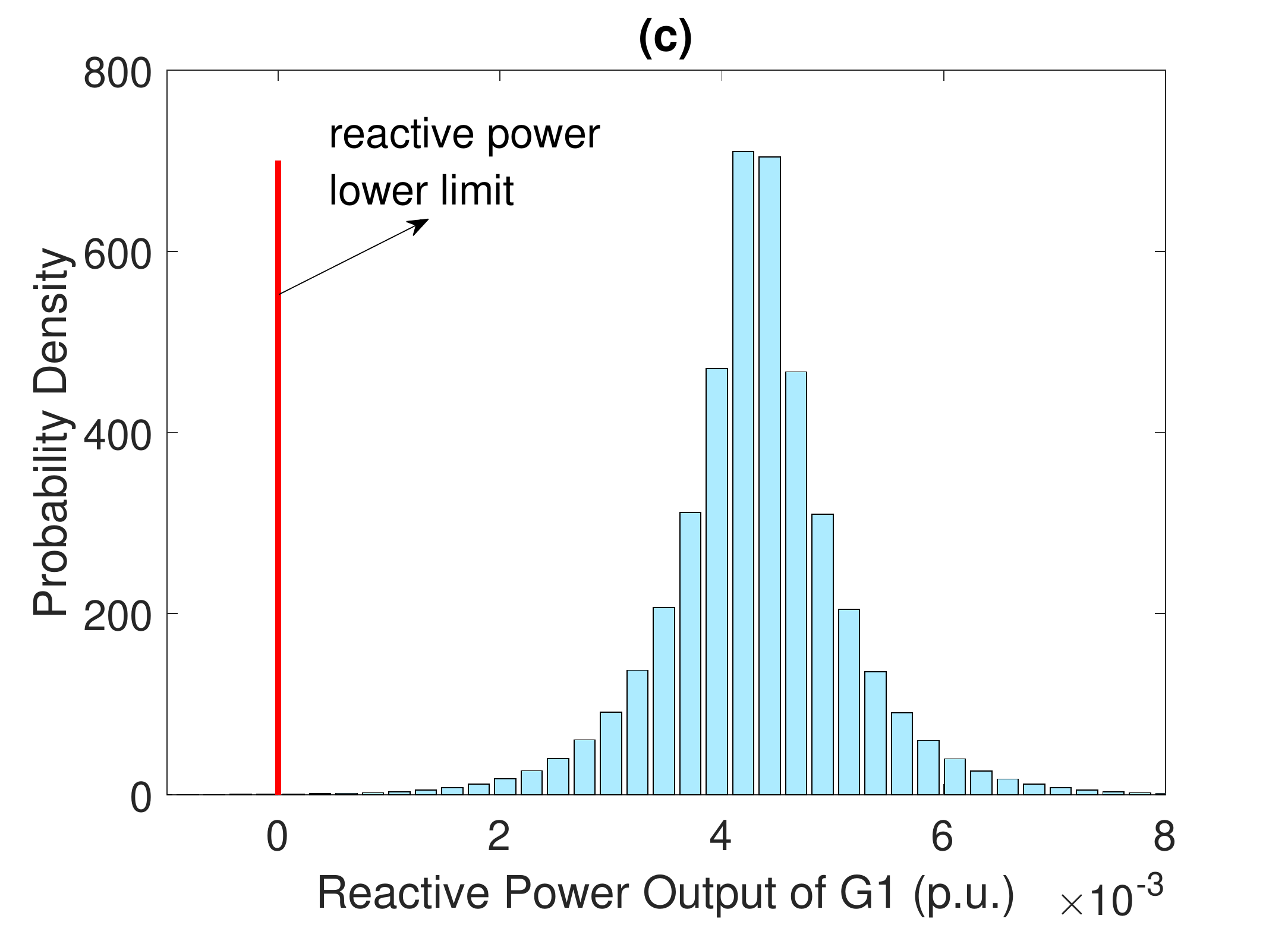}
		\includegraphics[width=1.7in]{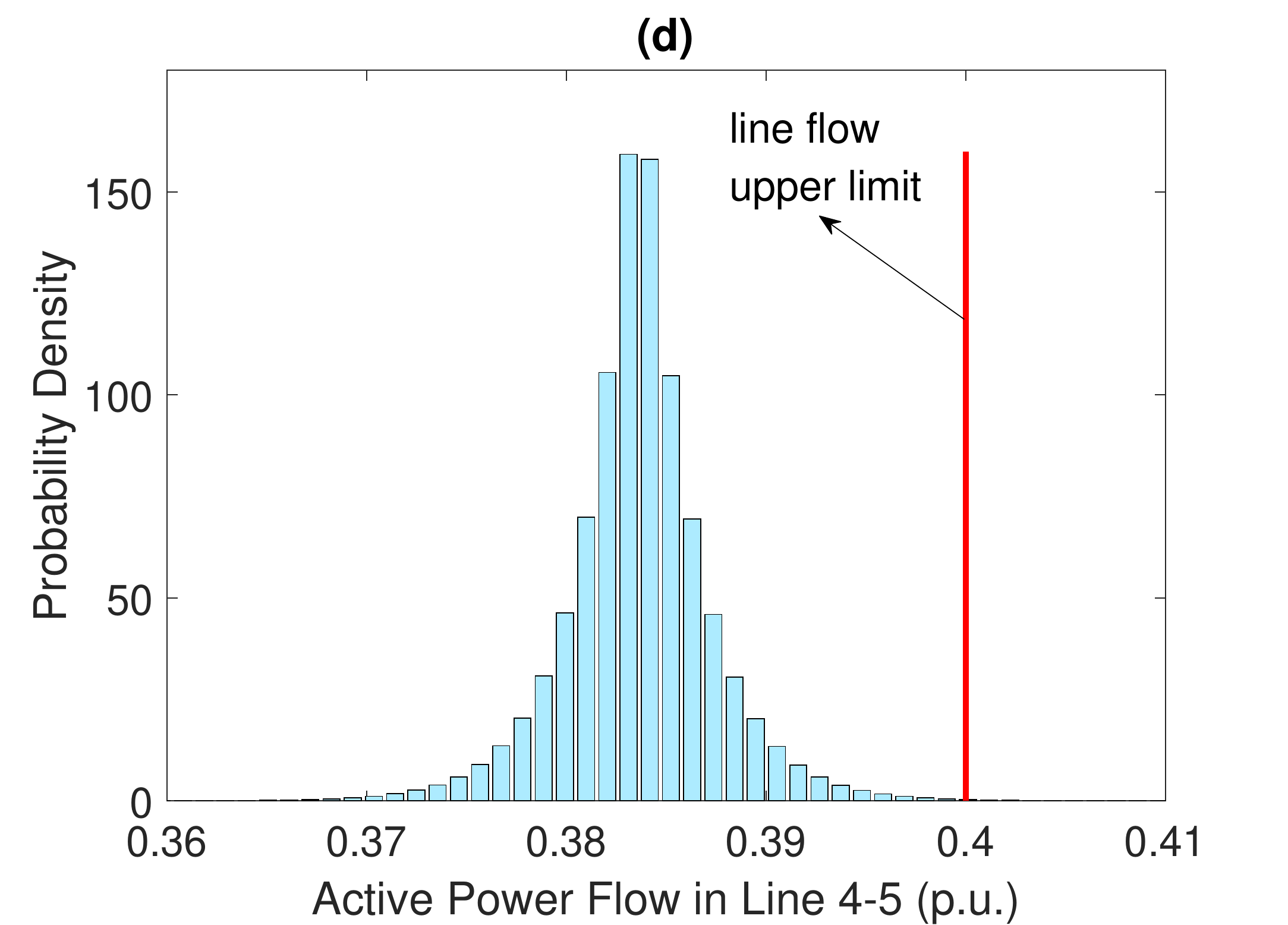}
		\caption{Normalized histogram of system reserve usage, voltage magnitude at bus 12, reactive power output of G1 and line flow in line 4-5 by MCS when operated at the strategy given by the proposed complete model} \label{plotdis}
	\end{figure}

\begin{figure}[!ht]\centering
		\includegraphics[width=1.7in]{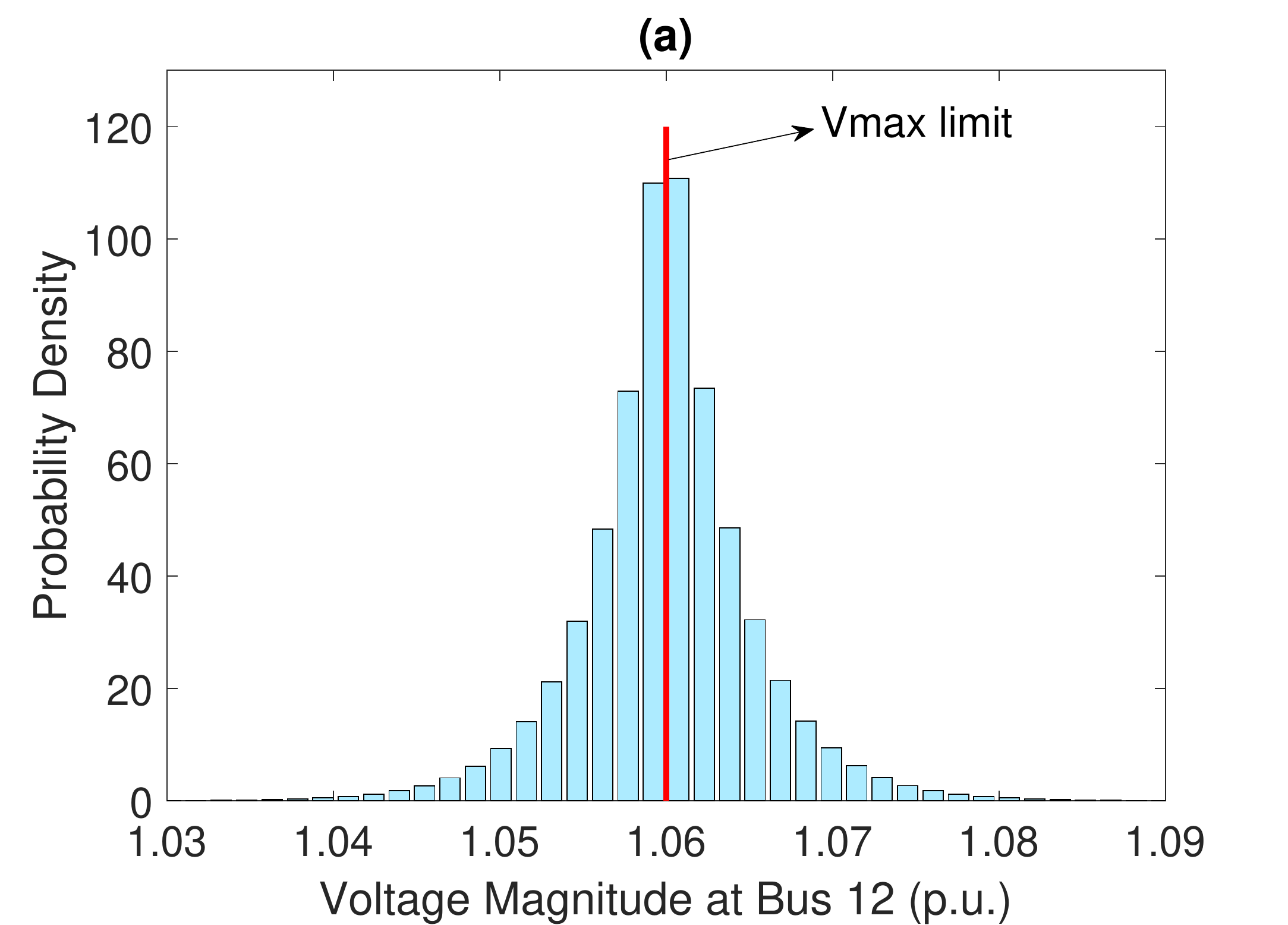}
		\includegraphics[width=1.7in]{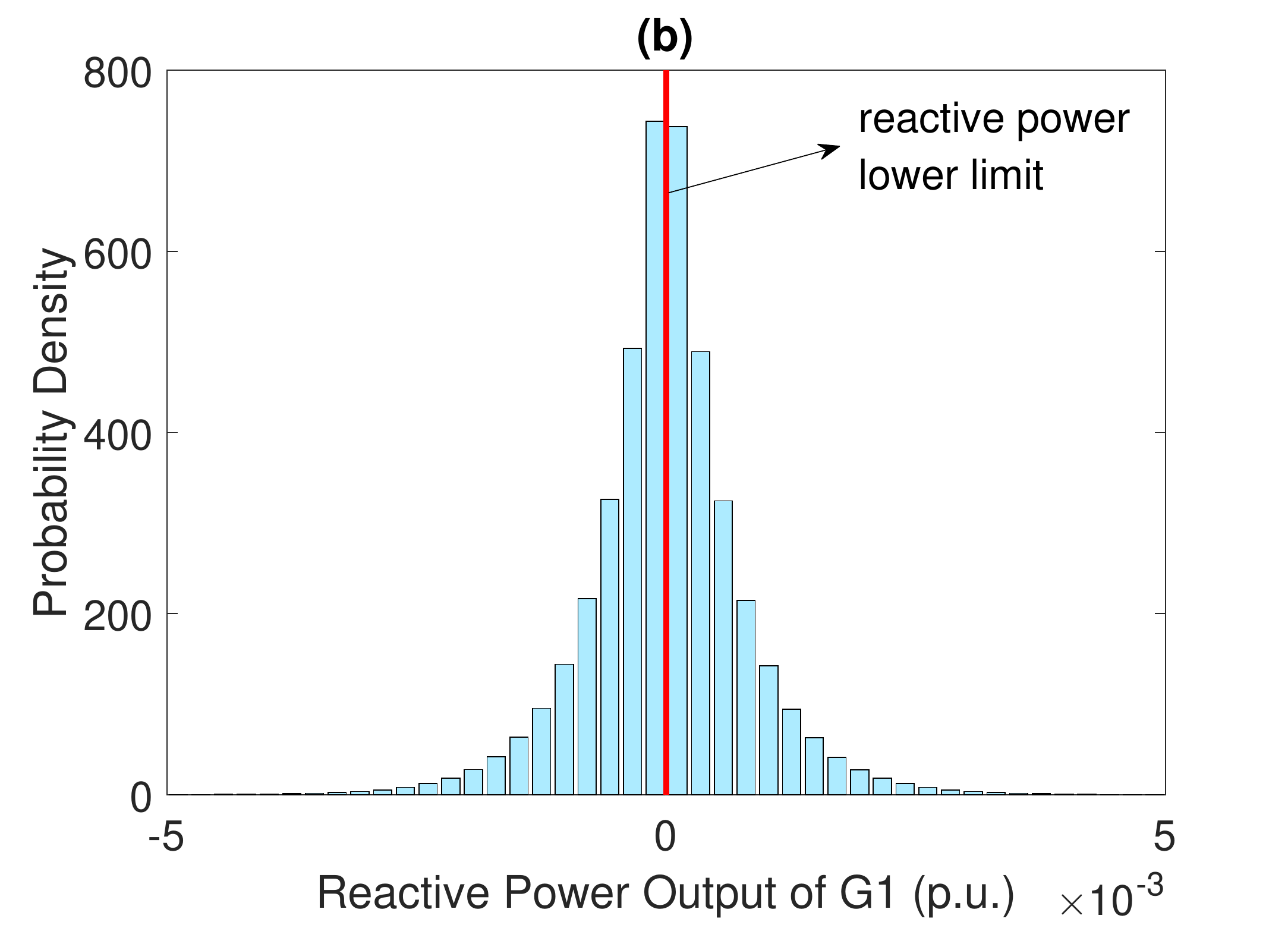}
		\caption{Normalized histogram of voltage magnitude at bus 12 (a) for the proposed model without voltage constraint and the normailized histogram of reactive power output at G1 (b) for the proposed model without reactive power constraint.} \label{plotdis_2}
	\end{figure}

	Table \ref{Strategy} presents the optimal operation strategy of the proposed model under different conditions with the comparison to the DC model. When operated with the strategy given by the complete model, all the chance constraints for reserve, voltage, reactive power and line flow are well guaranteed. For example, Fig. \ref{plotdis} shows the normalized histogram of system reserve usage, voltage magnitude at bus 12, reactive power output of G1 and line flow in line 4-5. It is shown those quantities safely distribute in the allowable ranges. If we exclude the voltage chance constraints (\ref{con6}) from the proposed model, the corresponding optimal generation strategy is also given in Table \ref{Strategy}. Compared with the strategy of complete model, the nominal active power and participation factors of generators do not have many differences, but the setting voltage magnitudes at generator buses increase significantly. The consequence of this change is shown in Fig. \ref{plotdis_2} (a). Under the voltage unconcerned strategy, the bus 12 will be exposed to the huge risk of over-voltage. On the other hand, if we drop the generator reactive power limits (\ref{conq}) from the complete model, the obtained operation strategy is also presented in Table \ref{Strategy}. Again, the nominal active power and participation factors only exhibit slight changes whereas the AVR setting points vary significantly. For example, the setting voltage of G1 decreases from 1.036 to 1.029. The consequence this change is illustrated in Fig. \ref{plotdis_2} (b) which shows G1 is under significant risk of crossing the under-excitation limit. Above discussion indicates that both voltage chance constraint (\ref{con6}) and the reactive power chance constraint (\ref{conq}) are essential for a safe and meaningful operation strategy under uncertainties. Table \ref{Strategy} further lists the strategy obtained by the DC model widely used in the literature \cite{bienstock2014chance,7332992,roald2015security,7478165,xie2016distributionally}. Except for the voltage and reactive power unawareness of DC model, more importantly, the nominal active power and participation factors exhibit non-negligible differences compared with those of the proposed model. This shows that a more accurate power flow model could significantly improve the operation strategy obtained by DC model.
	
	So it is important to ask how accurate is the approximate AC-OPF model in this paper? In our model (\ref{WDROPF_con}), we have an exact AC model (\ref{PF}) to govern the nominal operation point and the approximate LPF model (\ref{Vl})(\ref{Qg}) to calculate the incremental response under uncertainties. The accuracy, therefore, should lie between the LPF model (\ref{LPF}) and the exact AC model. The following group of tests is conducted under the operation strategy obtained by the proposed complete model shown in Table \ref{Strategy}. Table \ref{costcomp} shows the operation costs calculated using different models as the total VRE forecasting error varies from -32 MW to +32 MW. When the forecasting error is zero, our model coincides with the full-AC model, so the error of operation costs is zero. For other values of total forecasting error, the approximate AC model in this paper shows much higher accuracy than the LPF model in calculating the operation costs. Table \ref{voltcomp} and Table \ref{reactivecomp} further show the accuracy of calculating voltage magnitudes and reactive power outputs when the total VRE forecasting error is -32 MW. It is shown the accuracy of the approximate AC model in this paper is at least one-order-of-magnitude higher than that of the LPF model for calculating both voltage and reactive power. Since the approximate AC model in this paper inherits the formula (\ref{brflow}) of the LPF model, it has the same accuracy as the LPF model in calculating line MW flow but has higher accuracy than the DC model as shown in Table \ref{linecomp}. In summary, the proposed approximate AC model combines the accuracy of the full-AC model and the tractability of LPF model, which makes it a more attractive model for OPF under uncertainties.

	\begin{table*}[]
		\centering
		\caption{Optimal Operation Strategy of the Proposed Model under Different Conditions with Comparison to the DC Model}
		\label{Strategy}
		\begin{tabular}{c||ccc|ccc|ccc||ccc}
			\hline
			& \multicolumn{3}{c|}{the proposed complete model} & \multicolumn{3}{c|}{without voltage constratints} & \multicolumn{3}{c||}{without Qg constratints} & \multicolumn{3}{c}{DC Model} \\ \cline{2-13} 
			& Pg        & alpha     & Voltage     & Pg             & alpha          & Voltage         & Pg           & alpha        & Voltage        & Pg      & alpha   & Voltage  \\ \hline
			G1 & 0.562     & 0.061     & 1.036       & 0.561          & 0.061          & 1.058           & 0.566        & 0.066        & 1.029          & 0.555   & 0.062   & /        \\
			G2 & 0.339     & 0.292     & 1.026       & 0.339          & 0.293          & 1.048           & 0.343        & 0.309        & 1.019          & 0.338   & 0.195   & /        \\
			G3 & 0.538     & 0.000     & 1.024       & 0.539          & 0.000          & 1.045           & 0.546        & 0.000        & 1.017          & 0.454   & 0.000   & /        \\
			G4 & 0.376     & 0.535     & 1.032       & 0.376          & 0.535          & 1.055           & 0.370        & 0.525        & 1.033          & 0.386   & 0.548   & /        \\
			G5 & 0.079     & 0.112     & 0.995       & 0.078          & 0.111          & 1.018           & 0.070        & 0.100        & 0.986          & 0.137   & 0.195   & /        \\ \hline
		\end{tabular}
	\end{table*}

	\begin{table}[]
		\centering
		\caption{Operation Costs (\$) Calculated by Different Models}
		\label{costcomp}
		\begin{tabular}{c|c|c>{\bfseries}c|c>{\bfseries}c}
			\hline
			VRE forecast error & full-AC & this paper & error(\%) & LPF     & error(\%) \\ \hline
			-32                & 7502.07 & 7500.66    & 0.02\%    & 7440.90 & 0.82\%    \\
			-24                & 7180.80 & 7180.89    & 0.00\%    & 7121.23 & 0.83\%    \\
			-16                & 6863.44 & 6864.26    & -0.01\%   & 6804.70 & 0.86\%    \\
			-8                 & 6549.98 & 6550.76    & -0.01\%   & 6491.31 & 0.90\%    \\
			0                  & 6240.41 & 6240.41    & 0.00\%    & 6181.06 & 0.95\%    \\
			8                  & 5934.71 & 5933.19    & 0.03\%    & 5873.94 & 1.02\%    \\
			16                 & 5632.86 & 5629.11    & 0.07\%    & 5569.96 & 1.12\%    \\
			24                 & 5334.86 & 5328.17    & 0.13\%    & 5269.12 & 1.23\%    \\
			32                 & 5040.69 & 5030.36    & 0.21\%    & 4971.41 & 1.37\%    \\ \hline
		\end{tabular}
	\end{table}

	\begin{table}[]
		\centering
		\caption{Voltage Magnitudes (p.u.) Calculated by Different Models}
		\label{voltcomp}
		\begin{tabular}{c|c|c>{\bfseries}c|c>{\bfseries}c}
			\hline
			Bus No. & full-AC & this paper & error   & LPF  & error   \\ \hline
			4       & 1.019   & 1.019      & -0.0001 & 1.017 & -0.0023 \\
			5       & 1.025   & 1.025      & -0.0001 & 1.024 & -0.0012 \\
			7       & 1.002   & 1.002      & -0.0001 & 0.993 & -0.0094 \\
			9       & 0.998   & 0.998      & -0.0002 & 0.979 & -0.0192 \\
			10      & 0.997   & 0.997      & -0.0003 & 0.982 & -0.0153 \\
			11      & 1.015   & 1.014      & -0.0003 & 1.007 & -0.0073 \\
			12      & 1.028   & 1.028      & -0.0004 & 1.027 & -0.0012 \\
			13      & 1.018   & 1.017      & -0.0001 & 1.015 & -0.0025 \\
			14      & 0.995   & 0.995      & -0.0003 & 0.984 & -0.0114 \\ \hline
		\end{tabular}
	\end{table}

	\begin{table}[]
		\centering
		\caption{Reactive Power Output (p.u.) Calculated by Different Models}
		\label{reactivecomp}
		\begin{tabular}{c|c|c>{\bfseries}c|c>{\bfseries}c}
			\hline
			Unit No. & full-AC & this paper & error   & LPF   & error   \\ \hline
			G1       & 0.006   & 0.006      & -0.0006 & 0.005  & -0.0011 \\
			G2       & -0.101  & -0.103     & -0.0023 & -0.094 & 0.0067  \\
			G3       & 0.298   & 0.297      & -0.0005 & 0.301  & 0.0034  \\
			G4       & 0.211   & 0.213      & 0.0025  & 0.259  & 0.0480  \\
			G5       & -0.041  & -0.041     & -0.0001 & 0.011  & 0.0519  \\ \hline
		\end{tabular}
	\end{table}

	\begin{table}[]
		\centering
		\caption{Line MW Flow (p.u.) Calculated by Different Models}
		\label{linecomp}
		\begin{tabular}{c|c|c>{\bfseries}c|c>{\bfseries}c}
			\hline
			line & full-AC & this paper & error   & DC     & error   \\ \hline
			4-5   & -0.405  & -0.390     & 0.0149  & -0.340 & 0.0650  \\
			1-2   & 0.388   & 0.374      & -0.0139 & 0.323  & -0.0651 \\
			2-3   & 0.292   & 0.291      & -0.0019 & 0.261  & -0.0315 \\
			2-4   & 0.199   & 0.195      & -0.0035 & 0.148  & -0.0509 \\
			10-11 & -0.199  & -0.204     & -0.0058 & -0.185 & 0.0133  \\ \hline
		\end{tabular}
	\end{table}
	
	\subsection{IEEE 118-Bus System}
	On the IEEE 118-bus system, 18 wind farms are installed at bus 2, 5, 7, 13, 15, 21, 25, 28, 35, 45, 53, 58, 63, 75, 88, 95, 106 and 115, respectively. The capacity of each wind farm changes from 30 to 90 MW to test the methods under different levels of wind penetration and uncertainties. 
	
\begin{figure}[!ht]\centering
		\includegraphics[width=3.0in]{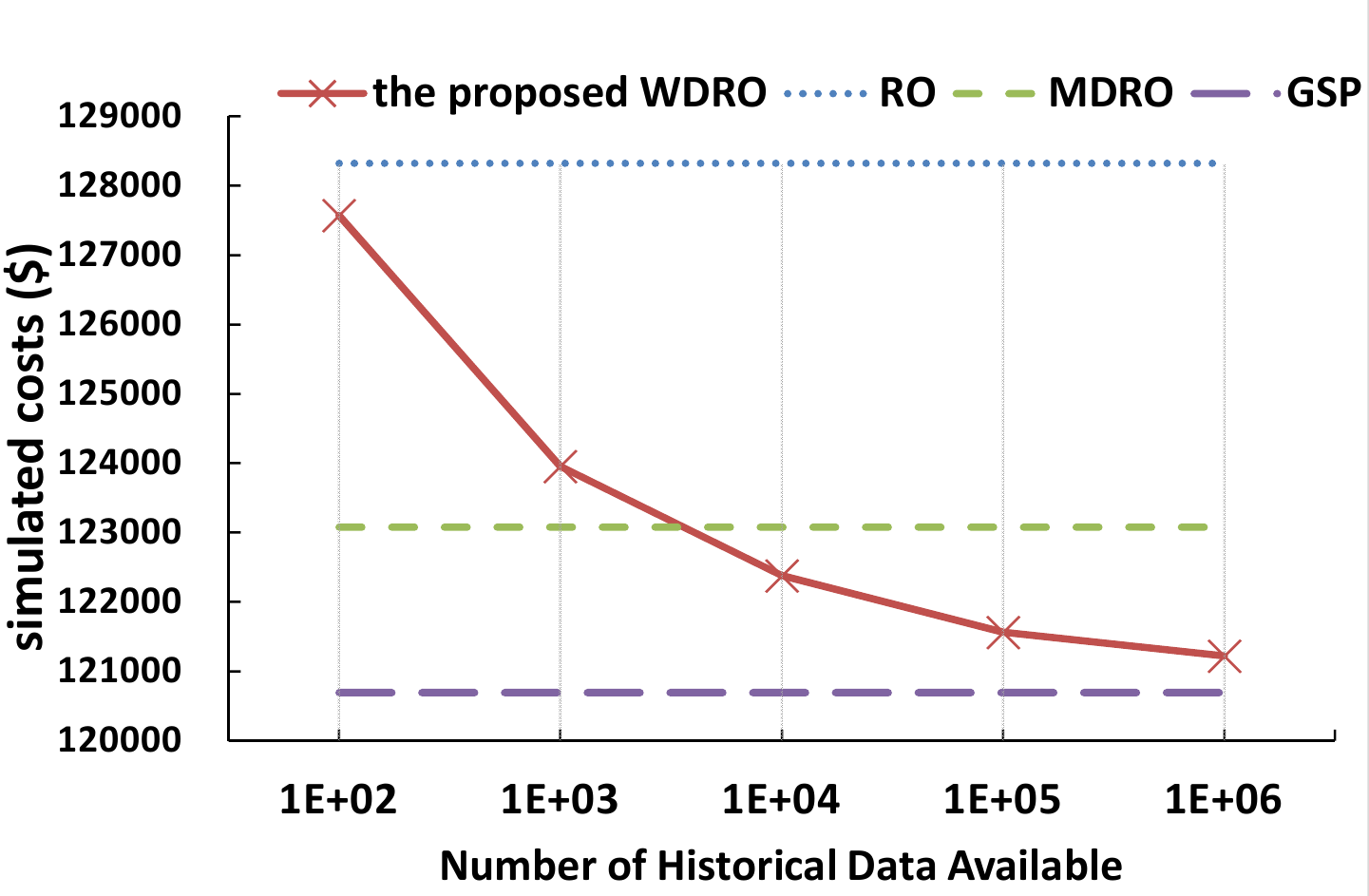}
		\caption{Evolution of simulated operational costs as more historical data is available with $30\times 18$ MW wind integration.} \label{datafigure1}
	\end{figure}

\begin{figure}[!ht]\centering
		\includegraphics[width=3.0in]{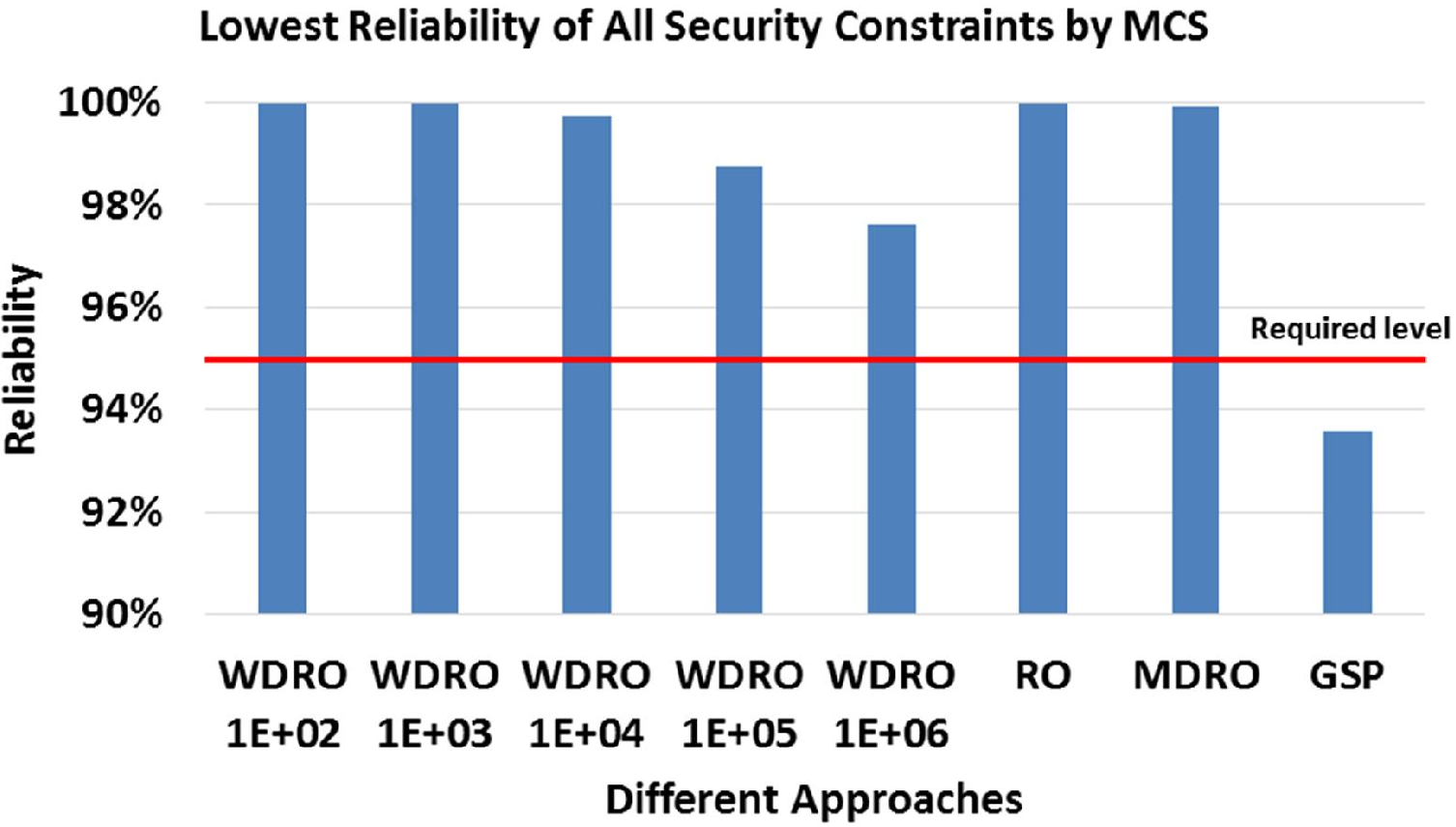}
		\caption{Comparison of the lowest reliability of all security constraints by Monte Carlo simulation} \label{datafigure2}
	\end{figure}

\begin{figure}[!ht]\centering
		\includegraphics[width=3.0in]{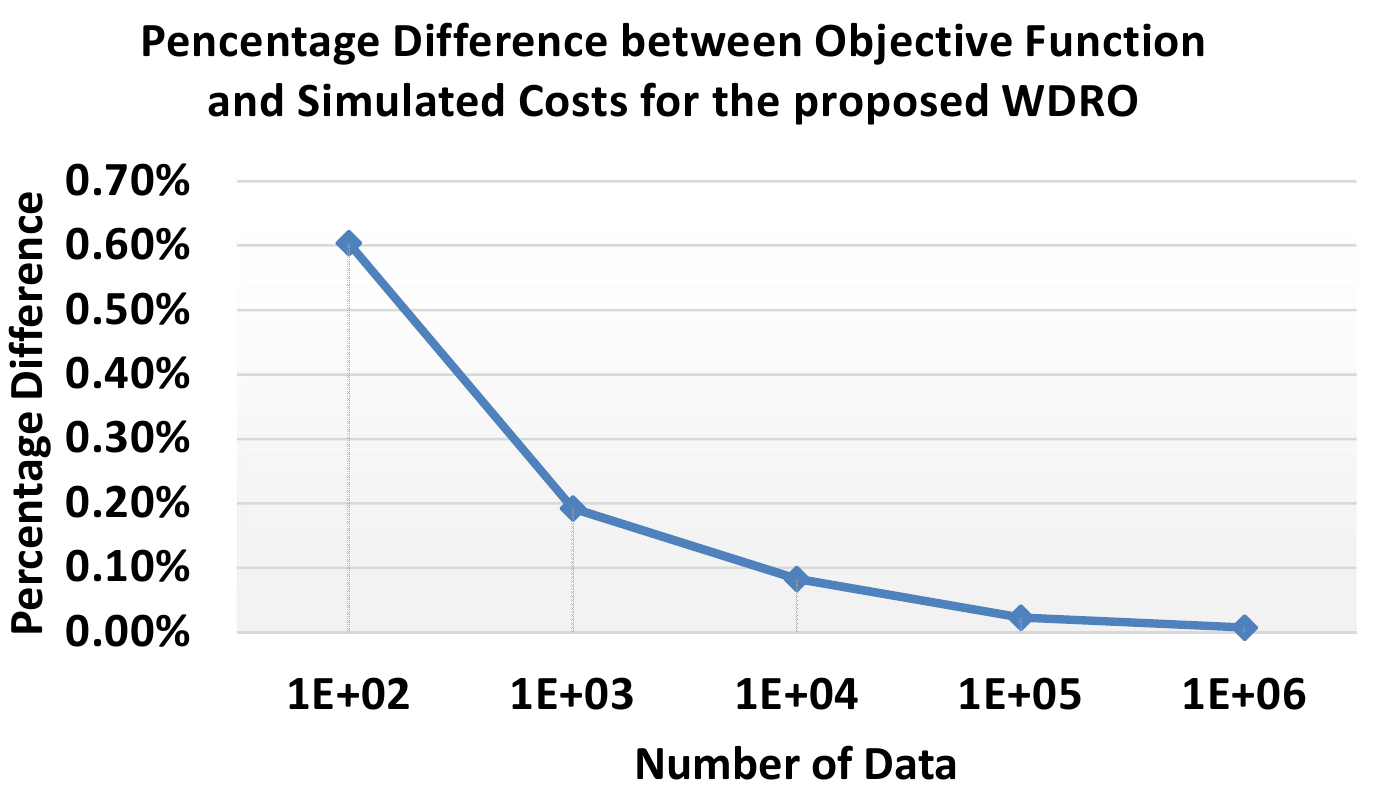}
		\caption{The procedure to determine the uncertainty set $\mathcal{U}$.} \label{datafigure3}
	\end{figure}

	In our implementation, the tolerable violation probability in (\ref{con5})$\sim$(\ref{con7}) is set to $\rho_1=\rho_2=\rho_3=\rho_4=0.05$, and the confidence level in (\ref{eps2}) is set to 0.9. The forecasting values of wind generation are set to be $50\%$ of the installed capacity. The available sample size ranges from $10^2$ to $10^6$ to showcase the data-exploiting feature. After solving each problem, Monte Carlo simulation with $10^7$ samples is employed to test the practical and out-of-sample performance. The major test results are summarized in Table \ref{my-label}. When the capacity of each wind farm is 30MW, all the tested methods obtain the corresponding optimal solutions. The evolution of the operational costs by Monte Carlo simulation is illustrated in Fig. \ref{datafigure1}. The costs of WDRO and MDRO lie between RO and GSP due to the fact that RO completely ignores the probabilistic information whereas the GSP assumes precise knowledge about the probability distribution. In other words, RO and GSP produce the most conservative and aggressive strategies, respectively. MDRO assumes partial knowledge, the first and second moments, of the probability distribution, which reduces conservatism compared with RO to some degree. In contrast, the proposed WDRO fully relies on the information told by the data at hand. When we are in short of data, it approaches the RO method to take a conservative decision. On the contrary, when the data is rich, it approaches the stochastic programming approach with complete information about the distribution. Fig. \ref{datafigure2} further compares the lowest reliability level of the security constraints. Although GSP obtains the strategy with lowest operation costs, it fails to guarantee the required reliability level. This is due to the deviation of the underlying true distribution from the Gauss assumption made in GSP method. All other methods, including RO, MDRO and the proposed WDRO ensure higher reliability level than required due to their ``robust'' nature. As more data is available, the proposed WDRO gradually and safely reduce the guaranteed reliability level to pursue higher economic efficiency. Note that the proposed WDRO considers the expected costs in the objective function (\ref{WDROPF_obj}) w.r.t. the worst-case distribution and the true distribution usually differs from the worst-case one, so the simulated costs are always upper bounded by the objective function as shown in Table \ref{my-label}. Fig. \ref{datafigure3} shows the evolution of percentage difference between the objective function and simulated costs for the proposed WDRO. As more data is available, the ambiguity set shrinks and worst-case distribution approaches the true distribution, and the difference between simulated costs and objective function diminishes. As shown in Table \ref{my-label}, when the installed capacity of wind power increases to $60\times 18$ MW, the RO along with the WDRO with $10^2$ and $10^3$ samples become infeasible. When the installed capacity increases to $90\times 18$ MW the MDRO and the WDRO with $10^4$ also become infeasible. In summary, we can order the methods according to their conservatism as: RO$>$WDRO($10^{2\sim 3}$)$>$MDRO$>$WDRO($10^{4\sim 6}$)$>$GSP.
	
	The computation of the proposed method consists of two stages. One is the preparation stage to construct the uncertainty set $\mathcal{U}$ in (\ref{RC}), and the other is the operation stage to solve the optimization problem (\ref{WDROPF_obj})(\ref{WDROPF_con}). In practice, the preparation stage can be done off-line and updated on weekly or monthly basis according to the new data available. In addition, the construction of the uncertainty set is highly parallelizable among different buses and transmission lines. The time on our 12-core workstation to construct uncertainty set for the whole 118-bus system is listed in Table \ref{my-label2}. The computation time grows linearly with the sample size $N$, and further speed-up is possible by leveraging distributed computation. The computation time for the on-line operation stage is listed in Table \ref{my-label}. In this stage, the solver time is nearly the same as the RO approach and less than that of MDRO and GSP methods. More importantly, the solver time of WDRO approach at the operation stage is irrelevant to the number of wind farms and amount of historical data at hand. Therefore, the proposed WDRO approach is highly scalable and efficient for real-time operation.
	
	\begin{table*}[]
		\centering
		\caption{Comparative Case Study on IEEE 118-bus System}
		\label{my-label}
		\scalebox{0.9}[0.9]{
			\begin{tabular}{c|l|ccccc|ccc}
				\hline
				\multirow{2}{*}{Wind Capacity} & \multicolumn{1}{c|}{\multirow{2}{*}{\diagbox{Performance}{Methods}}} & \multicolumn{5}{c|}{the proposed WDRO with different sample sizes}                                              & \multirow{2}{*}{RO}     & \multirow{2}{*}{MDRO}     & \multirow{2}{*}{GSP} \\
				& \multicolumn{1}{c|}{}                         & 1.E+02                      & 1.E+03                      & 1.E+04                      & 1.E+05     & 1.E+06     &                             &                             &                         \\ \hline
				\multirow{6}{*}{30*18 MW}      & objective (\$)                                     & 128345                       & 124192                       & 122483                       & 121588      & 121227      & 128355                       & 123168                       & 120839                   \\
				& simulated costs (\$)                               & 127574                       & 123953                       & 122381                       & 121560      & 121218      & 128322                       & 123081                       & 120694                   \\
				& up reverse (p.u.)                                   & 4.296                       & 2.504                       & 1.712                       & 1.302      & 1.130      & 4.672                      & 2.038                       & 0.828                   \\
				& down reverse (p.u.)                                   & 4.319                       & 2.490                       & 1.711                       & 1.299      & 1.129      & 4.690                       & 2.084                       & 0.906                   \\
				& reliability                                   & 99.99997\%                  & 99.98933\%                  & 99.74943\%                  & 98.75797\% & 97.61366\% & 100.00000\%                 & 99.93666\%                  & 93.56586\%              \\
				& cpu time (s)                                      & 0.31 & 0.27 & 0.31 & 0.23 & 0.29 & 0.45 & 0.27 & 0.29                    \\ \hline
				\multirow{6}{*}{60*18 MW}      & objective (\$)                                     & \multirow{6}{*}{Infeasible} & \multirow{6}{*}{Infeasible} & 115404                       & 113883      & 113017      & \multirow{6}{*}{Infeasible} & 116608                       & 112015                   \\
				& simulated costs (\$)                               &                             &                             & 115197                       & 113814      & 112990      &                             & 116692                       & 111867                   \\
				& up reverse (p.u.)                                     &                             &                             & 3.340                       & 2.672      & 2.258      &                             & 4.128                       & 1.657                   \\
				& down reverse (p.u.)                                   &                             &                             & 3.344                       & 2.677      & 2.259      &                             & 4.081                       & 1.738                   \\
				& reliability                                   &                             &                             & 99.70579\%                  & 98.90693\% & 97.61284\% &                             & 99.93669\%                  & 93.18895\%              \\
				& cpu time (s)                                      &                             &                             & 0.32                        & 0.44       & 0.27       &                             & 1.14                       & 0.26                    \\ \hline
				\multirow{6}{*}{90*18 MW}      & objective (\$)                                     & \multirow{6}{*}{Infeasible} & \multirow{6}{*}{Infeasible} & \multirow{6}{*}{Infeasible} & 106522      & 105369      & \multirow{6}{*}{Infeasible} & \multirow{6}{*}{Infeasible} & 103719                   \\
				& simulated cost (\$)                                 &                             &                             &                             & 106391      & 105315      &                             &                             & 103582                   \\
				& up reverse (p.u.)                                     &                             &                             &                             & 3.877      & 3.410      &                             &                             & 2.575                   \\
				& down reverse (p.u.)                                   &                             &                             &                             & 3.890      & 3.407      &                             &                             & 2.654                   \\
				& reliability                                   &                             &                             &                             & 98.72546\% & 97.66914\% &                             &                             & 93.75508\%              \\
				& cpu time (s)                                      &                             &                             &                             & 0.33       & 0.32       &                             &                             & 1.72                    \\ \hline
			\end{tabular}
		}
	\end{table*}

	\begin{table}[]
		\centering
		\caption{Time (s) to Construct Uncertainty Sets for 118-bus System}
		\label{my-label2}
		\begin{tabular}{|c|c|c|c|c|c|}
			\hline
			$N$        & 1.E+02 & 1.E+03 & 1.E+04 & 1.E+05 & 1.E+06 \\ \hline
			Time (s) & 0.67   & 0.82   & 1.64   & 16.48  & 291.64 \\ \hline
		\end{tabular}
	\end{table}

	\section{Conclusion}
	This paper has proposed a chance-constrained approximate AC-OPF under uncertainties based on distributionally robust optimization with Wasserstein metric. In order to overcome the flaws of the DC power flow model extensively used in stochastic OPF formulations, we have developed a tractable AC power flow model by integrating the full AC power flow and a linear power flow models. Numerical studies have demonstrated the proposed OPF formulation has improved precision as well as good numerical tractability. Without any assumption on the underlying probability distribution of the uncertainties, the proposed chance-constrained AC-OPF is completely driven by the available historical data. It extracts reliable probabilistic information from historical data to construct an ambiguity set of all possible distributions and immunizes the operation strategy against all distributions in the ambiguity set. The more data is available, the less conservative the solution is. In addition, special problem structures are properly exploited in the reformulation of the distributionally robust optimization problem to improve the scalability and efficiency of the numerical solution approach.

	\appendices
	
	\section*{Appendix}
	\begin{IEEEproof}[Proof of Lemma 2]
		Using the idea from section 5.2 of \cite{Esfahani2015Data}, for every $1\leq i \leq m$, we define
		\begin{equation}\label{defl1}
		\underline{l}_{i}(\tilde{\bm{\vartheta}})=
		\left\{
		\begin{aligned}
		&1 \quad &\text{if } \tilde{{\vartheta}}_{i}\leq -\sigma \\
		&-\infty  \quad &\ \text{otherwise } 
		\end{aligned}
		\right.
		\end{equation}
		\begin{equation}\label{defl2}
		\overline{l}_{i}(\tilde{\bm{\vartheta}})=
		\left\{
		\begin{aligned}
		&1  &\text{if } \tilde{{\vartheta}}_{i}\geq \sigma \\
		&-\infty   &\  \ \ \text{otherwise }. 
		\end{aligned}
		\right.
		\end{equation}
		Further define
		\begin{equation}\label{defl}
		l(\tilde{\bm{\vartheta}})=\text{max}\left\{ \underset{1\leq i \leq m}{\text{max}} \underline{l}_{i}(\tilde{\bm{\vartheta}}),\underset{1\leq i \leq m}{\text{max}} \overline{l}_{i}(\tilde{\bm{\vartheta}}),0  \right\}.
		\end{equation}
		It is easy to see that $l(\tilde{\bm{\vartheta}})$ is the indicator function of the compliment of $\mathcal{V}(\sigma)$, so $\mathbb{Q}[\tilde{\bm{\vartheta}} \notin \mathcal{V}(\sigma)]=\mathbb{E}_{\mathbb{Q}}[l(\tilde{\bm{\vartheta}})]$. In addition, $l(\tilde{\bm{\vartheta}})$ is the maximum of upper-semicontinuous functions, which implies $l(\tilde{\bm{\vartheta}})$ itself is upper-semicontinuous. Therefore, Lemma 1 applies:
		
		\begin{subequations}
			\begin{equation}
			\underset{\mathbb{Q}\in \hat{\mathcal{Q}}_{N}}{\text{sup}} \mathbb{Q}[\tilde{\bm{\vartheta}} \notin \mathcal{V}(\sigma)]
			= \underset{\mathbb{Q}\in \hat{\mathcal{Q}}_{N}}{\text{sup}} \mathbb{E}_{\mathbb{Q}}[l(\tilde{\bm{\vartheta}})] \qquad \ \
			\end{equation}    
			\begin{equation}\label{dc1}
			= \left\{
			\begin{array}{l}
			\underset{\lambda \geq 0,\bm{s} \in \mathbb{R}^{N}}{\text{inf}} 
			\ \lambda \cdot \epsilon+\frac{1}{N}\sum_{k=1}^{N}s_{k}\\
			\text{s.t.}\ \forall 1\leq k \leq N:\\
			\quad \ \underset{\bm{\vartheta} \in \Theta}{\text{sup}} \left( l({\bm{\vartheta}})-\lambda \norm{\bm{\vartheta}-\hat{\bm{\vartheta}}^{(k)}}_{1}  \right) \leq s_{k}
			\end{array}\right.
			\end{equation}    
			\begin{equation}\label{dc2}
			= \left\{
			\begin{array}{l}
			\underset{\lambda \geq 0,\bm{s} \geq \bm{0}}{\text{inf}} 
			\ \lambda \cdot \epsilon+\frac{1}{N}\sum_{k=1}^{N}s_{k}\\
			\text{s.t.}\ \forall 1\leq k \leq N, 1\leq i\leq m:\\
			\quad \ \underset{\bm{\vartheta} \in \Theta}{\text{sup}} \left( \underline{l}_{i}(\tilde{\bm{\vartheta}})-\lambda \norm{\bm{\vartheta}-\hat{\bm{\vartheta}}^{(k)}}_{1}  \right) \leq s_{k}\\
			\quad \ \underset{\bm{\vartheta} \in \Theta}{\text{sup}} \left( \overline{l}_{i}(\tilde{\bm{\vartheta}})-\lambda \norm{\bm{\vartheta}-\hat{\bm{\vartheta}}^{(k)}}_{1}  \right) \leq s_{k}
			\end{array}\right. 
			\end{equation}    
			\begin{equation}\label{dc3}
			= \left\{
			\begin{array}{l}
			\underset{\lambda \geq 0,\bm{s} \geq \bm{0}}{\text{inf}} 
			\ \lambda \cdot \epsilon+\frac{1}{N}\sum_{k=1}^{N}s_{k}\\
			\text{s.t.}\ \forall 1\leq k \leq N, 1\leq i\leq m:\qquad \qquad \quad
			\\
			\quad \ 1-\lambda(\hat{{\vartheta}}^{(k)}_{i}+\sigma)^{+} \leq s_{k}\\
			\quad \ 1-\lambda(\sigma-\hat{{\vartheta}}^{(k)}_{i})^{+} \leq s_{k}
			\end{array}\right.  
			\end{equation}    
			\begin{equation}\label{dc5}
			= \underset{\lambda \geq 0}{\text{inf}} \left\{
			\  \lambda \cdot \epsilon+\frac{1}{N}\sum_{k=1}^{N}\left(1-\lambda (\sigma-\norm{\hat{\bm{\vartheta}}^{({k})}}_{\infty})^{+}  \right)^{+} \right\} 
			\end{equation}
		\end{subequations}
		where (\ref{dc1}) is a direct application of Lemma 1; (\ref{dc2}) uses the the definition of $l(\tilde{\bm{\vartheta}})$ given in (\ref{defl}); (\ref{dc3}) explicitly evaluates the supremums in (\ref{dc2}) by exploiting the definitions of $\Theta$,  (\ref{defl1})(\ref{defl2}) and the $l_{1}$ norm; some direct calculation of (\ref{dc3}) leads to  (\ref{dc5}) which completes the proof.

	\end{IEEEproof}
	
	\vspace{-0.5cm}
	
	\bibliographystyle{ieeetr}
	
	\bibliography{DisRUC}
	% that's all folks

   \begin{IEEEbiographynophoto}
		{Chao Duan (S'14)}
		was born in Chongqing, China, in 1989. He received the B.S. degree in electrical engineering from Xi'an Jiaotong University, Xi'an, China,
		in 2012. He is currently pursuing the Ph.D. degree at Xi'an Jiaotong University,	Xi'an, China, and the University of Liverpool, Liverpool, U.K.
		His research interests are in stochastic optimization, stability analysis and robust control of power systems.
	\end{IEEEbiographynophoto}

	\begin{IEEEbiographynophoto}
		{Wanliang Fang}
		was born in Henan, China, in 1958. He received the B.S. and M.S. degrees from Xi'an Jiaotong University, Xi'an, China in 1982 and 1988, respectively, and the Ph.D. degree from Hong Kong Polytechnic University, HongKong, in 1999, all in electrical engineering.
		He is currently a Professor of electrical engineering at Xi'an Jiaotong University. His research interests include power system stability analysis and control, FACTS and HVDC.
	\end{IEEEbiographynophoto}

		\begin{IEEEbiographynophoto}
		{Lin Jiang (M'00)} received the B.Sc. and M.Sc. degrees from Huazhong University of Science and Technology (HUST), China, in 1992 and 1996; and the Ph.D. degree from the University of Liverpool, UK, in 2001, all in electrical engineering.
		He is currently a Reader in the University of Liverpool. His current research interests include control and analysis of power system, smart grid, and renewable energy.
	\end{IEEEbiographynophoto}

   \begin{IEEEbiographynophoto}
	{Li Yao}
	is currently pursuing the Ph.D. degree at Xi'an Jiaotong University,	Xi'an, China.
	His research interests are in distributionally robust optimization of power systems.
\end{IEEEbiographynophoto}

	\begin{IEEEbiographynophoto}
		{Jun Liu (S'09-M'10)}
		received the B.S. and Ph.D. degrees from Xi'an Jiaotong
		University, Xi'an, China, in 2004 and 2012, respectively, all in electrical engineering.
		He is currently an associate Professor in the Department of Electrical Engineering, Xi'an Jiaotong University. His research interests include renewable energy integration, power system operation and control, power system stability, HVDC, and FACTS.
	\end{IEEEbiographynophoto}

\end{document}